\documentclass[12pt,a4paper]{amsart}

\usepackage{lscape}
\usepackage{adjustbox}
\usepackage{hyperref}
\usepackage{ulem}

\newcommand {\fB}{{\mathfrak{B}}}
\newcommand {\fbgl}{{\mathfrak{bgl}}}
\newcommand {\fbr}{{\mathfrak{br}}}

\newcommand {\fbrj}{{\mathfrak{brj}}}
\newcommand {\fel}{{\mathfrak{el}}}

\newcommand {\fwk}{{\mathfrak{wk}}}

\usepackage{DLdef1}

\begin{document}

\title[
The roots of exceptional Lie superalgebras]{The roots of exceptional modular Lie superalgebras with Cartan matrix}

\author{Sofiane Bouarroudj${}^a$, 
Dimitry Leites${}^{a, b, *}$,\\ Alexander Lozhechnyk${}^{c}$, Jin Shang${}^a$}

\address{${}^a$New York University Abu Dhabi,
Division of Science and Mathematics, P.O. Box 129188, United Arab
Emirates; \{sofiane.bouarroudj, dl146, js8544\}@nyu.edu}
\address{${}^b$ Department of mathematics, Stockholm University, Roslagsv. 101: mleites@math.su.se\\
${}^{*}$ The corresponding author}
\address{${}^c$https://www.linkedin.com/
in/aleksandr-lozhechnik-40021263, alozhechhnik@gmail.com.}


\begin{abstract} For each of the exceptional Lie superalgebras with indecomposable Cartan matrix, we give the explicit list of its roots of and the corresponding Chevalley basis for one of the inequivalent Cartan matrices, the one corresponding to the greatest number of mutually orthogonal isotropic odd simple roots.

Our main tools: Grozman's Mathematica-based code SuperLie, and Python.
\end{abstract}

\thanks{S.B. and D.L. were partly supported by the grant AD 065 NYUAD. A.K was partly supported by Oberwolfach Leibniz Fellowship.}


\maketitle

\markboth{\itshape S.~Bouarroudj\textup{,} 
D.~ Leites\textup{,} A.~Lozhechnyk\textup{,} J.~Shang}{{\itshape Exceptional root systems}}

\thispagestyle{empty}

This paper is an Appendix to \cite{BGL} which contains classification of finite-dimensional Lie superalgebras  $\fg(A)$ with indecomposable Cartan matrix $A$ over an algebraically closed field $\Kee$ of characteristic $p>0$ together with all inequivalent Cartan matrices and the corresponding Dynkin--Kac graphs for each such Lie superalgebra.  The paper \cite{BGL}  contains  superdimensions of each Lie superalgebra  $\fg(A)$, a description of its even part $\fg(A)_\ev$, and of its odd part $\fg(A)_\od$ as a  $\fg(A)_\ev$-module. Recall that $A/a|B$ means that $A|B=\sdim \fg(A)$ and  
$a|B=\sdim \fg^{(1)}(A)/\fc$, where $\fg^{(1)}$ is the first derived and $\fc$ is the center.

Here, for each exceptional $\fg(A)$ we list its roots and root vectors  using Cartan matrix $A_m$ --- the one with the maximal number of not connected grey vertices in its Dynkin--Kac graph $D_m$. The set of simple roots related to $D_m$ is needed in relation with the Duflo--Serganova functor.

The tables below are analogs of descriptions of root systems in \cite{Bbk}. Recall that whereas each Lie algebra with indecomposable Cartan matrix over $\Cee$ (finite-dimensional or of polynomial growth, or even hyperbolic, see \cite{CCLL}) has just one Cartan matrix, each of the exceptional Lie \textbf{super}algebras  considered  has several Cartan matrices, see \cite{BGL,CCLL}; one Lie superalgebra we consider has 135 inequivalent Cartan matrices, dozens inequivalent Cartan matrices for one Lie superalgebra is usual. Therefore, to list sets of simple roots for all exceptional Lie superalgebras is hardly reasonable, so we list system of roots for each exceptional Lie superalgebra, but list \textit{simple roots} for just one Cartan matrix: for which the number of pairwise orthogonal isotropic simple roots is maximal (this Cartan matrix of $\fg$ is useful to define the defect of $\fg$). To find simple roots for the other Cartan matrices, one has to use odd reflections, see \cite{BGL}.

By $N\fg(A)$ we denote the realization of $\fg(A)$ corresponding to the $N$th Cartan matrix $A$ as listed in \cite{BGL}.  The odd root vectors are \fbox{boxed} and isotropic roots are \underline{underlined}.

\section{Roots and root vectors}

Let the $\pi_i$ be the fundamental weights relative a fixed system of simple roots. For any simple Lie algebra with Cartan matrix, let $R(\sum a_i\pi_i)$ denote both the irreducible representation  with highest weight $(a_1, \dots )$ and the respective module. The tautological module over the Lie algebra of series $\fsl$, $\fo$ or $\fsp$ is denoted by $\id:=R(\pi_1)$.

\ssec{$\fosp(4|2; a)$, $\fag(2)$,  and $\fab(3)$ for $p\geq 5$} The answer is the same as is well-known for $p=0$, namely (here $\fsl_i(2)$ is the $i$th copy of $\fsl(2)$ with the tautological module $\id_i$): 
\[
\begin{tabular}{|c|c|}
\hline
&\\[-10pt]
$\fosp(4|2; a)_\ev=\fsl_1(2)\oplus \fsl_2(2)\oplus\fsl_3(2)$&$\fosp(4|2; a)_\od=\id_1\boxtimes \id_2\boxtimes \id_3$\\[2pt]
$\fag(2)_\ev=\fsl(2)\oplus \fg(2)$&$\fag(2)_\od=\id\boxtimes R(\pi_1)$\\[2pt]
 $\fab(3)_\ev=\fsl(2)\oplus \fo(7)$&$\fab(3)_\od=\id\boxtimes R(\pi_1)$\\[2pt]
\hline
\end{tabular}
\]
Each of the other exceptional Lie superalgebras $\fg(A)$ with indecomposable Cartan matrix $A$ exists only in characteristics 2, 3 and 5. The Lie superalgebras $3\fg(2,3)$ and $1\fg(3,3)$ (indigenous to $p=3$) resemble $3\fag(2)$ and $6\fab(3)$ (existing for $p=0$ and any $p>3$), respectively, other exceptional Lie superalgebras $\fg(A)$ have no analogs except for two pairs $\fbrj(2;5)\leftrightarrow\fbrj(2;3)$ and
$\fel(5;5)\leftrightarrow\fel(5;3)$ which we consider one after the other for clarity.

\sssec{$1\fbrj(2;5)$  of $\sdim 10|12$}\label{brj25} We have $\fbrj(2; 5)_\ev = \fsp(4)\simeq\fbr(2; -1)$ (the Brown algebra)  and the $\fbrj(2; 5)_\ev$-module
$\fbrj(2; 5)_\od = R(\pi_1 + \pi_2)$  is irreducible. We consider the   Cartan matrix and basis elements
\[
\quad 
\begin{tabular}{|l|l|} 
\hline
&\\[-10pt]
 the root vectors&the roots\\
\hline
&\\[-8pt]
\fbox{$x_{1}$},\;\fbox{$x_{2}$}&\underline{$\alpha_1$},\; $\alpha_2$\\[3pt]
$x_3=\left[x_{1},\,x_{2}\right]$,\;
$x_4=\left[x_{2},\,x_{2}\right]$&$\alpha_1+\alpha_2,\ \ 2\alpha_2$\\[3pt]
$\fbox{$x_5$}=\left[x_{2},\,\left[x_{1},\,x_{2}\right]\right]$&$\alpha_1+2\alpha_2$\\[3pt]
$x_6=\left[\left[x_{1},\,x_{2}\right],\,\left[x_{2},\,x_{2}\right]
\right]$&$\alpha_1+3\alpha_2$\\[8pt]
$\begin{array}{@{}l@{}}
\fbox{$x_7$}=\left[\left[x_{1},\,x_{2}\right],\left[x_{2},\,\left[x_{1},\,x_{2}\right]\right]\right],\\
\fbox{$x_8$}=\left[\left[x_{2},\,x_{2}\right],\left[x_{2},\,\left[x_{1},\,x_{2}\right]\right]\right]
\end{array}
$&\underline{$2\alpha_1+3\alpha_2$},\ \ \underline{$\alpha_1+4\alpha_2$}\\[18pt]
$x_9=\left[\left[x_{1},\,x_{2}\right],\left[\left[x_{1},\,x_{2}\right],\,\left[x_{2},\,x_{2}\right]
 \right]\right]$&$2\alpha_1+4\alpha_2$\\[3pt]
$\fbox{$x_{10}$}=\left[\left[x_{2},\,\left[x_{1},\,x_{2}\right]
 \right],\,\left[\left[x_{1},\,x_{2}\right], \left[x_{2},\,x_{2}\right]\right]\right]$&\underline{$2\alpha_1+5\alpha_2$}\\[4pt]
 \hline
\end{tabular}
\]

\sssec{$1\fbrj(2;3)$ of $\sdim 10|8$}\label{brj23} We have $\fbrj(2; 3)_\ev = \fsp(4)\simeq\fbr(2;-1)$ (the Brown algebra) and the $\fbrj(2; 3)_\ev$-module
$\fbrj(2; 3)_\od = R(2\pi_2)$. We consider the   Cartan matrix and basis elements
\[
\quad 
\begin{tabular}{|l|l|} 
\hline
&\\[-10pt]
 the root vectors&the roots\\
\hline
&\\[-10pt]
\fbox{$x_{1}$},\;\fbox{$x_{2}$}&\underline{$\alpha_1$},\; $\alpha_2$\\[4pt]
$x_3=\left[x_{1},\,x_{2}
\right]$, \ $x_4=\left[x_{2},\,x_{2}\right]$&$\alpha_1+\alpha_2$, \ $2\alpha_2$\\[4pt]
$\fbox{$x_5$}=\left[x_{2},\,\left[x_{1},\,x_{2}
 \right]\right]$&$\alpha_1+2\alpha_2$\\[4pt]
$x_6=\left[\left[x_{1},\,x_{2}\right],\left[x_{2},\,x_{2}\right]\right]$&$\alpha_1+3\alpha_2$\\[4pt]
$\fbox{$x_7$}=\left[\left[x_{2},\,x_{2}\right],\left[x_{2},
\left[x_{1},\,x_{2}\right]\right]
 \right]$&\underline{$\alpha_1+4\alpha_2$}\\[4pt]
$x_8=\left[
\left[x_{1},\,x_{2}\right],\left[\left[x_{1},\,x_{2}\right],
\left[x_{2},\,x_{2}\right]\right] \right]$&$2\alpha_1+4\alpha_2$\\[4pt]
\hline
\end{tabular}
\]

\sssec{$1\fel(5;5)$ of $\sdim = 55|32$}\label{el(5;5)} We have  $\fel(5;5)_\ev = \fo(11)$ and $\fel(5;5)_\od =R(\pi_5)$ as $\fel(5;5)_\ev $-module. We consider the   Cartan matrix 
\begin{equation*} 
%
\end{equation*}
the Chevalley basis and the root system are:
\begin{equation*}
\label{basisEl55}
\footnotesize
\begin{tabularx}{\linewidth}{|>{\raggedright\let\newline\\\arraybackslash\hspace{0pt}}>{\hsize=1.4\hsize}X| >{\raggedright\let\newline\\\arraybackslash\hspace{0pt}}>{\hsize=0.6\hsize}X |}
\hline
&\\[-10pt]
 the root vectors&the roots\\
\hline
&\\[-10pt]
$x_1$,\ $x_2$,\ \fbox{$x_3$},\ \fbox{$x_4$},\ \fbox{$x_5$} & $\alpha_1$,\ $\alpha_2$,\ $\underline{\alpha_3}$,\ $\underline{\alpha_4}$,\ $\underline{\alpha_5}$ \\[2mm]
\hline
&\\[-10pt]
\fbox{$x_6$}=$[x_1,x_3]$,\ \fbox{$x_7$}=$[x_2,x_5]$,\ $x_8=[x_3,x_4]$, \ $x_9=[x_3,x_5]$,\ $x_{10}=[x_4,x_5]$ & $\underline{\alpha_1+\alpha_3}$,\ $\underline{\alpha_2+\alpha_5}$,\ $\alpha_3+\alpha_4$,\ $\alpha_3+\alpha_5$,\ $\alpha_4+\alpha_5$ \\[2mm]
\hline
&\\[-10pt]
$x_{11}=[x_3,[x_2,x_5]]$,\ $x_{12}=[x_4,[x_1,x_3]]$,\ $x_{13}=[x_4,[x_2,x_5]]$,\ $x_{14}=[x_5,[x_1,x_3]]$,\ \fbox{$x_{15}$}=$[x_5,[x_3,x_4]]$ & $\alpha_2+\alpha_3+\alpha_5$,\ $\alpha_1+\alpha_3+\alpha_4$,\ $\alpha_2+\alpha_4+\alpha_5$,\ $\alpha_1+\alpha_3+\alpha_5$,\ $\underline{\alpha_3+\alpha_4+\alpha_5}$ \\
\hline
\end{tabularx}
\end{equation*}

\sssec{$7\fel(5;3)$ of $\sdim = 39|32$}\label{el(5;3)} We have  $\fel(5;3)_\ev = \fo(9)\oplus \fsl(2)$ and $\fel(5;3)_\od = R(\pi_4)\boxtimes \id$ as $\fel(5;3)_\ev$-module.
We consider the   Cartan matrix and Chevalley basis elements
\begin{equation*}

\end{equation*}
\begin{equation*}
\begin{tabularx}{\linewidth}{|>{\raggedright\let\newline\\\arraybackslash\hspace{0pt}}>{\hsize=1.39\hsize}X| >{\raggedright\let\newline\\\arraybackslash\hspace{0pt}}>{\hsize=0.61\hsize}X |}
\hline
&\\[-9pt]
the root vectors&the roots\\
\hline
&\\[-9pt]
$\boxed{x_{1}}$,
$x_{2}$,
$\boxed{x_{3}}$,
$\boxed{x_{4}}$,
$\boxed{x_{5}}$,
&
\underline{$\alpha_{1}$},
$\alpha_{2}$,
\underline{$\alpha_{3}$},
\underline{$\alpha_{4}$},
\underline{$\alpha_{5}$},
\\[5pt]
\hline
&\\[-9pt]
$x_{6}=[x_{1}, x_{4}]$,\newline
$\boxed{x_{7} }=[x_{2}, x_{4}]$,\newline
$x_{8}=[x_{3}, x_{4}]$,\newline
$x_{9}=[x_{3}, x_{5}]$,
&
$\alpha_{1} + \alpha_{4}$,\newline
\underline{$\alpha_{2} + \alpha_{4}$},\newline
$\alpha_{3} + \alpha_{4}$,\newline
$\alpha_{3} + \alpha_{5}$,
\\[5pt]
\hline
&\\[-9pt]
$x_{10}=[x_{2}, [x_{1}, x_{4}]]$,\newline
$\boxed{x_{11} }=[x_{3}, [x_{1}, x_{4}]]$,\newline
$x_{12}=[x_{3}, [x_{2}, x_{4}]]$,\newline
$\boxed{x_{13} }=[x_{5}, [x_{3}, x_{4}]]$,
&
$\alpha_{1} + \alpha_{2} + \alpha_{4}$,\newline
\underline{$\alpha_{1} + \alpha_{3} + \alpha_{4}$},\newline
$\alpha_{2} + \alpha_{3} + \alpha_{4}$,\newline
\underline{$\alpha_{3} + \alpha_{4} + \alpha_{5}$},
\\
\hline
\end{tabularx}
\end{equation*}

\sssec{$1\fg(1,6)$  of $\sdim 21|14$}\label{g(1,6)} We have $\fg(1,6)_\ev = \fsp(6)$ and $\fg(1,6)_\od = R(\pi_3)$ as $\fg(1,6)_\ev$-module. We consider the Cartan matrix $ %
$ and the Chevalley basis elements 
\iffalse }\\
\text{even $|$ odd }\end{matrix} \; \tiny
\begin{array}{l}
x_{1},\;|\;x_{2},\;x_{3}\\
x_5=\left[x_{2},\,x_{2}\right],\;|\; x_4=
\left[x_{1},\,x_{2}\right],\;
\\
 x_6=\left[x_{2},\,x_{3}\right],\;
x_7=\left[x_{2},\left[x_{1},\,x_{2}\right]\right]\;|\;\\
 x_8=\left[x_{3},\left[x_{1},\,x_{2}\right]\right],\;|\; x_9=\left[x_{3},\left[x_{2},\,x_{2}\right]\right]\\
x_{10}=\left[\left[x_{1},\,x_{2}\right],\left[x_{1},\,x_{2}\right]\right],\;|\;
x_{11}=\left[\left[x_{1},\,x_{2}\right],\left[x_{2},\,x_{3}\right]\right]\\
\;|\;x_{12}=\left[\left[x_{1},\,x_{2}\right],\left[x_{3},
\left[x_{1},\,x_{2}\right]\right]
 \right]\\
 x_{13}=\left[
\left[x_{2},\,x_{3}\right],\left[x_{2},
\left[x_{1},\,x_{2}\right]\right]
 \right]\;|\;\\
 x_{14}=\left[
\left[x_{2}, \left[x_{1},\,x_{2}\right]\right]
 ,\,\left[x_{3}, \left[x_{1},\,x_{2}\right]\right]
 \right]\;|\;\\
 \;|\;x_{15}=\left[
\left[x_{3}, \left[x_{2},\,x_{2}\right]\right]
 ,\,\left[\left[x_{1},\,x_{2}
\right], \left[x_{1},\,x_{2}\right]\right]
 \right]\\
 x_{16}=\left[
\left[x_{3}, \left[x_{2},\,x_{2}\right]\right]
 ,\,\left[\left[x_{1},\,x_{2}
\right], \left[x_{3},\left[x_{1},\,x_{2}\right]
\right]\right]\right]\;|\;
\end{array}
\end{equation}\fi

\begin{equation*}
\begin{tabularx}{\linewidth}{|>{\raggedright\let\newline\\\arraybackslash\hspace{0pt}}>{\hsize=1.4\hsize}X| >{\raggedright\let\newline\\\arraybackslash\hspace{0pt}}>{\hsize=0.6\hsize}X |}
\hline
&\\[-10pt]
the root vectors&the roots\\[3pt]
\hline
&\\[-10pt]
$x_{1}$,
$\boxed{x_{2}}$,
$\boxed{x_{3}}$,
&
$\alpha_{1}$,
$\alpha_{2}$,
\underline{$\alpha_{3}$},
\\[4pt]
\hline
&\\[-10pt]
$\boxed{x_{4} }=[x_{1}, x_{2}]$,
$x_{5}=[x_{2}, x_{2}]$,
$x_{6}=[x_{2}, x_{3}]$,
&
$\alpha_{1} + \alpha_{2}$,
$2  \alpha_{2}$,
$\alpha_{2} + \alpha_{3}$,
\\[4pt]
\hline
&\\[-10pt]
$x_{7}=[x_{2}, [x_{1}, x_{2}]]$,
$x_{8}=[x_{3}, [x_{1}, x_{2}]]$,
$\boxed{x_{9} }=[x_{3}, [x_{2}, x_{2}]]$,
&
$2  \alpha_{2} + \alpha_{1}$,
$\alpha_{1} + \alpha_{2} + \alpha_{3}$,
\underline{$2  \alpha_{2} + \alpha_{3}$},
\\[4pt]
\hline
&\\[-10pt]
$x_{10}=[[x_{1}, x_{2}], [x_{1}, x_{2}]]$,
$\boxed{x_{11} }=[[x_{1}, x_{2}], [x_{2}, x_{3}]]$,
&
$2  \alpha_{1} + 2  \alpha_{2}$,
$2  \alpha_{2} + \alpha_{1} + \alpha_{3}$,
\\[4pt]
\hline
&\\[-10pt]
$\boxed{x_{12} }=[[x_{1}, x_{2}], [x_{3}, [x_{1}, x_{2}]]]$,
$x_{13}=[[x_{2}, x_{3}], [x_{2}, [x_{1}, x_{2}]]]$,
&
\underline{$2  \alpha_{1} + 2  \alpha_{2} + \alpha_{3}$},
$3  \alpha_{2} + \alpha_{1} + \alpha_{3}$,
\\[4pt]
\hline
&\\[-10pt]
$x_{14}=[[x_{2}, [x_{1}, x_{2}]], [x_{3}, [x_{1}, x_{2}]]]$,
&
$2  \alpha_{1} + 3  \alpha_{2} + \alpha_{3}$,
\\[4pt]
\hline
&\\[-10pt]
$\boxed{x_{15} }=[[x_{3}, [x_{2}, x_{2}]], [[x_{1}, x_{2}], [x_{1}, x_{2}]]]$,
&
\underline{$2  \alpha_{1} + 4  \alpha_{2} + \alpha_{3}$},
\\[4pt]
\hline
&\\[-10pt]
$x_{16}=[[[x_{1}, x_{2}], [x_{2}, x_{3}]], [[x_{1}, x_{2}], [x_{2}, x_{3}]]]$
&
$2  \alpha_{1} + 2  \alpha_{3} + 4  \alpha_{2}$
\\[4pt]
\hline
\end{tabularx}
\end{equation*}

\sssec{$2\fg(2,3)$ of $\sdim 12/10|14$}\label{SSSg33} We have $\fg(2,3)_\ev = \fgl(3) \oplus \fsl(2)$ and $\fg(2,3)_\od = \fpsl(3)\boxtimes\id$  as $\fg(2,3)_\ev$-module. Clearly, $(\fg^{(1)}(2,3)/\fc)_\ev = \fpsl(3) \oplus \fsl(2)$. We consider the Cartan matrix
$%
$ and the Chevalley basis elements
\begin{equation*}
\begin{tabularx}{\linewidth}{|>{\raggedright\let\newline\\\arraybackslash\hspace{0pt}}>{\hsize=1.4\hsize}X| >{\raggedright\let\newline\\\arraybackslash\hspace{0pt}}>{\hsize=0.6\hsize}X |}
\hline
&\\[-10pt]
the root vectors&the roots\\[4pt]
\hline
&\\[-10pt]
$\boxed{x_{1}}$,
$\boxed{x_{2}}$,
$\boxed{x_{3}}$,
&
\underline{$\alpha_{1}$},
\underline{$\alpha_{2}$},
\underline{$\alpha_{3}$},
\\[4pt]
\hline
&\\[-10pt]
$x_{4}=[x_{1}, x_{3}]$,
$x_{5}=[x_{2}, x_{3}]$,
&
$\alpha_{1} + \alpha_{3}$,
$\alpha_{2} + \alpha_{3}$,
\\[4pt]
\hline
&\\[-10pt]
$\boxed{x_{6} }=[x_{2}, [x_{1}, x_{3}]]$,
&
$\alpha_{1} + \alpha_{2} + \alpha_{3}$,
\\[4pt]
\hline
&\\[-10pt]
$x_{7}=[[x_{1}, x_{3}], [x_{2}, x_{3}]]$,
&
$2  \alpha_{3} + \alpha_{1} + \alpha_{2}$,
\\[4pt]
\hline
&\\[-10pt]
$\boxed{x_{8} }=[[x_{1}, x_{3}], [x_{2}, [x_{1}, x_{3}]]]$,
$\boxed{x_{9} }=[[x_{2}, x_{3}], [x_{2}, [x_{1}, x_{3}]]]$,
&
\underline{$2  \alpha_{1} + 2  \alpha_{3} + \alpha_{2}$},
\underline{$2  \alpha_{2} + 2  \alpha_{3} + \alpha_{1}$},
\\[4pt]
\hline
&\\[-10pt]
$x_{10}=[[x_{2}, [x_{1}, x_{3}]], [x_{2}, [x_{1}, x_{3}]]]$,
&
$2  \alpha_{1} + 2  \alpha_{2} + 2  \alpha_{3}$,
\\[4pt]
\hline
&\\[-10pt]
$\boxed{x_{11} }=[[x_{2}, [x_{1}, x_{3}]], [[x_{1}, x_{3}], [x_{2}, x_{3}]]]$
&
\underline{$2  \alpha_{1} + 2  \alpha_{2} + 3  \alpha_{3}$}
\\[4pt]
\hline
\end{tabularx}
\end{equation*}

\sssec{$2\fg(2,6)$  of $\sdim 36/34|20$}\label{g26} We have $\fg(2,6)_\ev = \fgl(6)$ and $\fg(2,6)_\od = R(\pi_3)$ as $\fg(2,6)_\ev$-module. Clearly, $(\fg^{(1)}(2,6)/\fc)_\ev = \fpsl(6)$.
We consider the Cartan matrix 
\begin{equation*}
%
\end{equation*}
and the Chevalley basis elements
\begin{equation*}
\begin{tabularx}{\linewidth}{|>{\raggedright\let\newline\\\arraybackslash\hspace{0pt}}>{\hsize=1.2\hsize}X| >{\raggedright\let\newline\\\arraybackslash\hspace{0pt}}>{\hsize=0.8\hsize}X |}
\hline
&\\[-10pt]
the root vectors&the roots\\[3pt]
\hline
&\\[-10pt]
$x_{1}$,
$\boxed{x_{2}}$,
$\boxed{x_{3}}$,
$\boxed{x_{4}}$,
$x_{5}$,
&
$\alpha_{1}$,
\underline{$\alpha_{2}$},
\underline{$\alpha_{3}$},
\underline{$\alpha_{4}$},
$\alpha_{5}$,
\\[3pt]
\hline
&\\[-10pt]
$\boxed{x_{6} }=[x_{1}, x_{2}]$,
$x_{7}=[x_{2}, x_{3}]$,
$x_{8}=[x_{2}, x_{4}]$,
$x_{9}=[x_{3}, x_{4}]$,
$\boxed{x_{10} }=[x_{3}, x_{5}]$,
&
\underline{$\alpha_{1} + \alpha_{2}$},
$\alpha_{2} + \alpha_{3}$,
$\alpha_{2} + \alpha_{4}$,
$\alpha_{3} + \alpha_{4}$,
\underline{$\alpha_{3} + \alpha_{5}$},
\\[3pt]
\hline
&\\[-10pt]
$x_{11}=[x_{3}, [x_{1}, x_{2}]]$,
$x_{12}=[x_{4}, [x_{1}, x_{2}]]$,
$\boxed{x_{13} }=[x_{4}, [x_{2}, x_{3}]]$,
$x_{14}=[x_{5}, [x_{2}, x_{3}]]$,
$x_{15}=[x_{5}, [x_{3}, x_{4}]]$,
&
$\alpha_{1} + \alpha_{2} + \alpha_{3}$,
$\alpha_{1} + \alpha_{2} + \alpha_{4}$,
\underline{$\alpha_{2} + \alpha_{3} + \alpha_{4}$},
$\alpha_{2} + \alpha_{3} + \alpha_{5}$,
$\alpha_{3} + \alpha_{4} + \alpha_{5}$,
\\
\hline
\end{tabularx}
\end{equation*}

\sssec{7$\fg(3,3)$  of $\sdim 23/21|16$}\label{SSSg33} 
Let $\spin_7:=R(\pi_3)$. We have 
$\fg(3,3)_\ev=(\fo(7)\oplus \Kee z)\oplus \Kee d$ and 
$\fg(3,3)_\od=(\spin_7)_+\oplus (\spin_7)_-$ as $\fg(3,3)_\ev$-module. The action of $d$ --- the outer derivative of $\fg(3,3)^{(1)}$ --- separates the  identical $\fo(7)$-modules $\spin_7$  by acting on these modules as the scalar multiplication by $\pm1$, as indicated by subscripts, $z$ spans the center of $\fg(3,3)$. We consider the Cartan matrix $%
$
and the Chevalley basis elements
\begin{equation*}
\begin{tabularx}{\linewidth}{|>{\raggedright\let\newline\\\arraybackslash\hspace{0pt}}>{\hsize=1.4\hsize}X| >{\raggedright\let\newline\\\arraybackslash\hspace{0pt}}>{\hsize=0.6\hsize}X |}
\hline
&\\[-10pt]
the root vectors&the roots\\
\hline
&\\[-10pt]
$\boxed{x_{1}}$,
$\boxed{x_{2}}$,
$x_{3}$,
$\boxed{x_{4}}$,
&
\underline{$\alpha_{1}$},
\underline{$\alpha_{2}$},
$\alpha_{3}$,
\underline{$\alpha_{4}$},
\\[3pt]
\hline
&\\[-10pt]
$x_{5}=[x_{1}, x_{2}]$,
$\boxed{x_{6} }=[x_{2}, x_{3}]$,
$x_{7}=[x_{2}, x_{4}]$,
&
$\alpha_{1} + \alpha_{2}$,
\underline{$\alpha_{2} + \alpha_{3}$},
$\alpha_{2} + \alpha_{4}$,
\\[3pt]
\hline
&\\[-10pt]
$x_{8}=[x_{3}, [x_{1}, x_{2}]]$,
$\boxed{x_{9} }=[x_{4}, [x_{1}, x_{2}]]$,
$x_{10}=[x_{4}, [x_{2}, x_{3}]]$,
&
$\alpha_{1} + \alpha_{2} + \alpha_{3}$,
\underline{$\alpha_{1} + \alpha_{2} + \alpha_{4}$},
$\alpha_{2} + \alpha_{3} + \alpha_{4}$,
\\[3pt]
\hline
&\\[-10pt]
$\boxed{x_{11} }=[x_{4}, [x_{3}, [x_{1}, x_{2}]]]$,
$\boxed{x_{12} }=[[x_{1}, x_{2}], [x_{2}, x_{3}]]$,
&
\underline{$\alpha_{1} + \alpha_{2} + \alpha_{3} + \alpha_{4}$},
\underline{$2  \alpha_{2} + \alpha_{1} + \alpha_{3}$},
\\[3pt]
\hline
&\\[-10pt]
$x_{13}=[[x_{1}, x_{2}], [x_{3}, [x_{1}, x_{2}]]]$,
$x_{14}=[[x_{2}, x_{4}], [x_{3}, [x_{1}, x_{2}]]]$,
&
$2  \alpha_{1} + 2  \alpha_{2} + \alpha_{3}$,
$2  \alpha_{2} + \alpha_{1} + \alpha_{3} + \alpha_{4}$,
\\[3pt]
\hline
&\\[-10pt]
$\boxed{x_{15} }=[[x_{3}, [x_{1}, x_{2}]], [x_{4}, [x_{1}, x_{2}]]]$,
&
\underline{$2  \alpha_{1} + 2  \alpha_{2} + \alpha_{3} + \alpha_{4}$},
\\[3pt]
\hline
&\\[-10pt]
$x_{16}=[[x_{4}, [x_{1}, x_{2}]], [[x_{1}, x_{2}], [x_{2}, x_{3}]]]$,
&
$2  \alpha_{1} + 3  \alpha_{2} + \alpha_{3} + \alpha_{4}$,
\\[3pt]
\hline
&\\[-10pt]
$x_{17}=[[x_{4}, [x_{3}, [x_{1}, x_{2}]]], [[x_{1}, x_{2}], [x_{2}, x_{3}]]]$
&
$2  \alpha_{1} + 2  \alpha_{3} + 3  \alpha_{2} + \alpha_{4}$
\\[3pt]
\hline
\end{tabularx}
\end{equation*}

\sssec{$2\fg(3,6)$ of $\sdim 36|40$} We have $\fg(3,6)_\ev = \fsp(8)$ and $\fg(3,6)_\od = R(\pi_3)$ s $\fg(3,6)_\ev$-module. We consider the Cartan matrix 
\begin{equation*}
%
\end{equation*}
and the Chevalley basis elements
\begin{equation*}
\footnotesize
\begin{tabularx}{\linewidth}{|>{\raggedright\let\newline\\\arraybackslash\hspace{0pt}}>{\hsize=1.5\hsize}X| >{\raggedright\let\newline\\\arraybackslash\hspace{0pt}}>{\hsize=0.5\hsize}X |}
\hline
&\\[-9pt]
the root vectors&the roots\\[5pt]
\hline
&\\[-9pt]
$\boxed{x_{1}}$,
$\boxed{x_{2}}$,
$\boxed{x_{3}}$,
$\boxed{x_{4}}$,
&
\underline{$\alpha_{1}$},
\underline{$\alpha_{2}$},
$\alpha_{3}$,
\underline{$\alpha_{4}$},
\\[5pt]
\hline
&\\[-9pt]
$x_{5}=[x_{1}, x_{2}]$,\newline
$x_{6}=[x_{2}, x_{3}]$,\newline
$x_{7}=[x_{3}, x_{3}]$,\newline
$x_{8}=[x_{3}, x_{4}]$,
&
$\alpha_{1} + \alpha_{2}$,\newline
$\alpha_{2} + \alpha_{3}$,\newline
$2  \alpha_{3}$,\newline
$\alpha_{3} + \alpha_{4}$,
\\[5pt]
\hline
&\\[-9pt]
$\boxed{x_{9} }=[x_{3}, [x_{1}, x_{2}]]$,\newline
$\boxed{x_{10} }=[x_{3}, [x_{2}, x_{3}]]$,\newline
$\boxed{x_{11} }=[x_{4}, [x_{2}, x_{3}]]$,\newline
$\boxed{x_{12} }=[x_{4}, [x_{3}, x_{3}]]$,
&
$\alpha_{1} + \alpha_{2} + \alpha_{3}$,\newline
\underline{$2  \alpha_{3} + \alpha_{2}$},\newline
\underline{$\alpha_{2} + \alpha_{3} + \alpha_{4}$},\newline
\underline{$2  \alpha_{3} + \alpha_{4}$},
\\
\hline
\end{tabularx}
\end{equation*}

\ssec{$6\fg(4,3)$ of $\sdim 24|26$}\label{Lg4,3}  We have $\fg(4,3)_\ev = \fsp(6) \oplus \fsl(2)$ and $\fg(4,3)_\od = R(\pi_2) \boxtimes \id$ as $\fg(4,3)_\ev$-module. We consider the Cartan matrix
 \begin{equation*}
%
\\[15pt]
\end{equation*}
 and the Chevalley basis elements 
 \begin{equation*}
\begin{tabularx}{\linewidth}{|>{\raggedright\let\newline\\\arraybackslash\hspace{0pt}}>{\hsize=1.4\hsize}X| >{\raggedright\let\newline\\\arraybackslash\hspace{0pt}}>{\hsize=0.6\hsize}X |}
\hline
&\\[-10pt]
the root vectors&the roots\\[4pt]
\hline
&\\[-10pt]
$\boxed{x_{1}}$,\
$\boxed{x_{2}}$,\
$\boxed{x_{3}}$,\
$\boxed{x_{4}}$,
&
\underline{$\alpha_{1}$},\
\underline{$\alpha_{2}$},\
\underline{$\alpha_{3}$},\
\underline{$\alpha_{4}$},
\\[4pt]
\hline
&\\[-10pt]
$x_{5}=[x_{1}, x_{2}]$,\newline 
$x_{6}=[x_{2}, x_{3}]$,\newline 
$x_{7}=[x_{3}, x_{4}]$,
&
$\alpha_{1} + \alpha_{2}$,\newline 
$\alpha_{2} + \alpha_{3}$,\newline 
$\alpha_{3} + \alpha_{4}$,
\\[4pt]
\hline
&\\[-10pt]
$\boxed{x_{8} }=[x_{3}, [x_{1}, x_{2}]]$,\newline
$x_{9}=[x_{4}, [x_{2}, x_{3}]]$,
&
\underline{$\alpha_{1} + \alpha_{2} + \alpha_{3}$},\newline 
$\alpha_{2} + \alpha_{3} + \alpha_{4}$,
\\
\hline
\end{tabularx}
\end{equation*}

\sssec{$13\fg(8,3)$ of $\sdim 55|50$}\label{e83} We have $\fg(8, 3)_\ev = \ff(4) \oplus \fsl(2)$ and $\fg(8, 3)_\od = R(\pi_4) \boxtimes \id$ as $\fg(8, 3)_\ev$-module.
 We consider the Cartan matrix
\begin{equation*} 
%
\end{equation*}
 and the Chevalley basis elements
\begin{equation*}
\footnotesize
\begin{tabularx}{\linewidth}{|>{\raggedright\let\newline\\\arraybackslash\hspace{0pt}}>{\hsize=1.4\hsize}X| >{\raggedright\let\newline\\\arraybackslash\hspace{0pt}}>{\hsize=0.6\hsize}X |}
\hline
&\\[-10pt]
the root vectors&the roots\\
\hline
&\\[-10pt]
$x_{1}$,
$\boxed{x_{2}}$,
$\boxed{x_{3}}$,
$\boxed{x_{4}}$,
$\boxed{x_{5}}$,
&
$\alpha_{1}$,
\underline{$\alpha_{2}$},
\underline{$\alpha_{3}$},
\underline{$\alpha_{4}$},
\underline{$\alpha_{5}$},
\\[4pt]
\hline
&\\[-10pt]
$\boxed{x_{6} }=[x_{1}, x_{2}]$,
$\boxed{x_{7} }=[x_{1}, x_{3}]$,
$x_{8}=[x_{2}, x_{3}]$,\newline
$x_{9}=[x_{2}, x_{4}]$,
$x_{10}=[x_{4}, x_{5}]$,
&
\underline{$\alpha_{1} + \alpha_{2}$},
\underline{$\alpha_{1} + \alpha_{3}$},
$\alpha_{2} + \alpha_{3}$,
$\alpha_{2} + \alpha_{4}$,
$\alpha_{4} + \alpha_{5}$,
\\[4pt]
\hline
&\\[-10pt]
$x_{11}=[x_{3}, [x_{1}, x_{2}]]$,
$x_{12}=[x_{4}, [x_{1}, x_{2}]]$,\newline 
$\boxed{x_{13} }=[x_{4}, [x_{2}, x_{3}]]$,
$\boxed{x_{14} }=[x_{5}, [x_{2}, x_{4}]]$,
&
$\alpha_{1} + \alpha_{2} + \alpha_{3}$,
$\alpha_{1} + \alpha_{2} + \alpha_{4}$,
\underline{$\alpha_{2} + \alpha_{3} + \alpha_{4}$},
\underline{$\alpha_{2} + \alpha_{4} + \alpha_{5}$},
\\
\hline
\end{tabularx}
\end{equation*}

\sssec{$2\fg(4,6)$ of $\sdim 66|32$}\label{g46} We have $\fg(4,6)_\ev = \fo(12)$ and $\fg(4,6)_\od = R(\pi_5)$ as $\fg(4,6)_\ev$-module. We consider the Cartan matrix
 \begin{equation*} 
%
\end{equation*}
 and the Chevalley basis elements
\begin{equation*}
\footnotesize
\begin{tabularx}{\linewidth}{|>{\raggedright\let\newline\\\arraybackslash\hspace{0pt}}>{\hsize=1.35\hsize}X| >{\raggedright\let\newline\\\arraybackslash\hspace{0pt}}>{\hsize=0.65\hsize}X |}
\hline
&\\[-10pt]
the root vectors&the roots\\[4pt]
\hline
$x_{1}$,
$x_{2}$,
$\boxed{x_{3}}$,
$\boxed{x_{4}}$,
$\boxed{x_{5}}$,
$x_{6}$,
&
$\alpha_{1}$,
$\alpha_{2}$,
\underline{$\alpha_{3}$},
\underline{$\alpha_{4}$},
\underline{$\alpha_{5}$},
$\alpha_{6}$,
\\[4pt]
\hline
&\\[-10pt]
$x_{7}=[x_{1}, x_{2}]$,
$\boxed{x_{8} }=[x_{2}, x_{3}]$,
$x_{9}=[x_{3}, x_{4}]$,\newline
$x_{10}=[x_{3}, x_{5}]$,
$x_{11}=[x_{4}, x_{5}]$,
$\boxed{x_{12} }=[x_{4}, x_{6}]$,
&
$\alpha_{1} + \alpha_{2}$,
\underline{$\alpha_{2} + \alpha_{3}$},
$\alpha_{3} + \alpha_{4}$,
$\alpha_{3} + \alpha_{5}$,
$\alpha_{4} + \alpha_{5}$,
\underline{$\alpha_{4} + \alpha_{6}$},
\\[4pt]
\hline
&\\[-10pt]
$\boxed{x_{13} }=[x_{3}, [x_{1}, x_{2}]]$,
$x_{14}=[x_{4}, [x_{2}, x_{3}]]$,\newline
$x_{15}=[x_{5}, [x_{2}, x_{3}]]$,
$\boxed{x_{16} }=[x_{5}, [x_{3}, x_{4}]]$,\newline
$x_{17}=[x_{6}, [x_{3}, x_{4}]]$,
$x_{18}=[x_{6}, [x_{4}, x_{5}]]$,
&
\underline{$\alpha_{1} + \alpha_{2} + \alpha_{3}$},
$\alpha_{2} + \alpha_{3} + \alpha_{4}$,
$\alpha_{2} + \alpha_{3} + \alpha_{5}$,
\underline{$\alpha_{3} + \alpha_{4} + \alpha_{5}$},
$\alpha_{3} + \alpha_{4} + \alpha_{6}$,
$\alpha_{4} + \alpha_{5} + \alpha_{6}$,
\\[4pt]
\hline
&\\
\hline
\end{tabularx}
\end{equation*}

\sssec{$4\fg(6,6)$  of $\sdim 78|64$}\label{g66} We have $\fg(6,6)_\ev= \fo(13)$ and $\fg(6,6)_\od = \spin_{13}:=R(\pi_6)$ as $\fg(6,6)_\ev$-module. We consider the Cartan matrix 
 \begin{equation*}
%
\end{equation*}
and the Chevalley basis elements
\begin{equation*}
\small
\begin{tabularx}{\linewidth}{|>{\raggedright\let\newline\\\arraybackslash\hspace{0pt}}>{\hsize=1.27\hsize}X| >{\raggedright\let\newline\\\arraybackslash\hspace{0pt}}>{\hsize=0.73\hsize}X |}
\hline
&\\[-10pt]
the root vectors&the roots\\[4pt]
\hline
&\\[-10pt]
$\boxed{x_{1}}$,
$\boxed{x_{2}}$,
$\boxed{x_{3}}$,
$\boxed{x_{4}}$,
$\boxed{x_{5}}$,
$x_{6}$,
&
\underline{$\alpha_{1}$},
\underline{$\alpha_{2}$},
\underline{$\alpha_{3}$},
\underline{$\alpha_{4}$},
\underline{$\alpha_{5}$},
$\alpha_{6}$,
\\[4pt]
\hline
&\\[-10pt]
$x_{7}=[x_{1}, x_{2}]$,
$x_{8}=[x_{2}, x_{3}]$,
$x_{9}=[x_{3}, x_{4}]$,\newline
$x_{10}=[x_{3}, x_{5}]$,
$x_{11}=[x_{4}, x_{5}]$,
$\boxed{x_{12} }=[x_{4}, x_{6}]$,
&
$\alpha_{1} + \alpha_{2}$,
$\alpha_{2} + \alpha_{3}$,
$\alpha_{3} + \alpha_{4}$,
$\alpha_{3} + \alpha_{5}$,
$\alpha_{4} + \alpha_{5}$,
\underline{$\alpha_{4} + \alpha_{6}$},
\\[4pt]
\hline
&\\[-10pt]
$\boxed{x_{13} }=[x_{3}, [x_{1}, x_{2}]]$,
$\boxed{x_{14} }=[x_{4}, [x_{2}, x_{3}]]$,\newline
$\boxed{x_{15} }=[x_{5}, [x_{2}, x_{3}]]$,
$\boxed{x_{16} }=[x_{5}, [x_{3}, x_{4}]]$,\newline
$x_{17}=[x_{6}, [x_{3}, x_{4}]]$,
$x_{18}=[x_{6}, [x_{4}, x_{5}]]$,
&
\underline{$\alpha_{1} + \alpha_{2} + \alpha_{3}$},
\underline{$\alpha_{2} + \alpha_{3} + \alpha_{4}$},
\underline{$\alpha_{2} + \alpha_{3} + \alpha_{5}$},
\underline{$\alpha_{3} + \alpha_{4} + \alpha_{5}$},
$\alpha_{3} + \alpha_{4} + \alpha_{6}$,
$\alpha_{4} + \alpha_{5} + \alpha_{6}$,
\\
\hline
\end{tabularx}
\end{equation*}

\begin{equation*}
\small
\begin{tabularx}{\linewidth}{|>{\raggedright\let\newline\\\arraybackslash\hspace{0pt}}>{\hsize=1.27\hsize}X| >{\raggedright\let\newline\\\arraybackslash\hspace{0pt}}>{\hsize=0.73\hsize}X |}
\hline
&\\[-10pt]
$\boxed{x_{47} }=[[[x_{3}, x_{5}], [x_{4}, [x_{2}, x_{3}]]], [[x_{4}, x_{6}], [x_{3}, [x_{1}, x_{2}]]]]$,
$\boxed{x_{48} }=[[[x_{3}, x_{5}], [x_{4}, [x_{2}, x_{3}]]], [[x_{4}, x_{6}], [x_{5}, [x_{2}, x_{3}]]]]$,
$\boxed{x_{49} }=[[[x_{4}, x_{5}], [x_{3}, [x_{1}, x_{2}]]], [[x_{4}, x_{6}], [x_{5}, [x_{2}, x_{3}]]]]$,
&
\underline{$2  \alpha_{2} + 2  \alpha_{4} + 3  \alpha_{3} + \alpha_{1} + \alpha_{5} + \alpha_{6}$},
\underline{$2  \alpha_{2} + 2  \alpha_{4} + 2  \alpha_{5} + 3  \alpha_{3} + \alpha_{6}$},
\underline{$2  \alpha_{2} + 2  \alpha_{3} + 2  \alpha_{4} + 2  \alpha_{5} + \alpha_{1} + \alpha_{6}$},
\\[4pt]
\hline
&\\[-10pt]
$x_{50}=[[[x_{4}, x_{6}], [x_{3}, [x_{1}, x_{2}]]], [[x_{4}, [x_{2}, x_{3}]], [x_{5}, [x_{2}, x_{3}]]]]$,
$x_{51}=[[[x_{4}, x_{6}], [x_{5}, [x_{2}, x_{3}]]], [[x_{3}, [x_{1}, x_{2}]], [x_{5}, [x_{3}, x_{4}]]]]$,
$x_{52}=[[[x_{4}, x_{6}], [x_{5}, [x_{3}, x_{4}]]], [[x_{4}, [x_{2}, x_{3}]], [x_{5}, [x_{2}, x_{3}]]]]$,
&
$2  \alpha_{4} + 3  \alpha_{2} + 3  \alpha_{3} + \alpha_{1} + \alpha_{5} + \alpha_{6}$,
$2  \alpha_{2} + 2  \alpha_{4} + 2  \alpha_{5} + 3  \alpha_{3} + \alpha_{1} + \alpha_{6}$,
$2  \alpha_{2} + 2  \alpha_{5} + 3  \alpha_{3} + 3  \alpha_{4} + \alpha_{6}$,
\\[4pt]
\hline
&\\
\hline
\end{tabularx}
\end{equation*}

\sssec{$5\fg(8,6)$  of $\sdim 133|56$}\label{g86} We have $\fg(8,6)_\ev = \fe(7)$ and $\fg(8,6)_\od = R(\pi_1)$ as $\fg(8,6)_\ev$-module. We consider the following Cartan matrix and the Chevalley basis elements
\begin{equation*}

\end{equation*}
\begin{landscape}
\begin{tabularx}{\linewidth}{|>{\raggedright\let\newline\\\arraybackslash\hspace{0pt}}>{\hsize=1.35\hsize}X| >{\raggedright\let\newline\\\arraybackslash\hspace{0pt}}>{\hsize=0.65\hsize}X |}
\hline
&\\[-10pt]
the root vectors&the roots\\[4pt]
\hline
$x_{1}$,
$\boxed{x_{2}}$,
$\boxed{x_{3}}$,
$\boxed{x_{4}}$,
$x_{5}$,
$x_{6}$,
$x_{7}$,
&
$\alpha_{1}$,
\underline{$\alpha_{2}$},
\underline{$\alpha_{3}$},
\underline{$\alpha_{4}$},
$\alpha_{5}$,
$\alpha_{6}$,
$\alpha_{7}$,
\\[4pt]
\hline
&\\[-10pt]
$\boxed{x_{8} }=[x_{1}, x_{3}]$,
$x_{9}=[x_{2}, x_{3}]$,
$x_{10}=[x_{2}, x_{4}]$,
$x_{11}=[x_{3}, x_{4}]$,
$\boxed{x_{12} }=[x_{4}, x_{5}]$,
$x_{13}=[x_{5}, x_{6}]$,
$x_{14}=[x_{6}, x_{7}]$,
&
\underline{$\alpha_{1} + \alpha_{3}$},
$\alpha_{2} + \alpha_{3}$,
$\alpha_{2} + \alpha_{4}$,
$\alpha_{3} + \alpha_{4}$,
\underline{$\alpha_{4} + \alpha_{5}$},
$\alpha_{5} + \alpha_{6}$,
$\alpha_{6} + \alpha_{7}$,
\\[4pt]
\hline
&\\[-10pt]
$x_{15}=[x_{2}, [x_{1}, x_{3}]]$,
$x_{16}=[x_{4}, [x_{1}, x_{3}]]$,
$\boxed{x_{17} }=[x_{4}, [x_{2}, x_{3}]]$,
$x_{18}=[x_{5}, [x_{2}, x_{4}]]$,
$x_{19}=[x_{5}, [x_{3}, x_{4}]]$,
$\boxed{x_{20} }=[x_{6}, [x_{4}, x_{5}]]$,
$x_{21}=[x_{7}, [x_{5}, x_{6}]]$,
&
$\alpha_{1} + \alpha_{2} + \alpha_{3}$,
$\alpha_{1} + \alpha_{3} + \alpha_{4}$,
\underline{$\alpha_{2} + \alpha_{3} + \alpha_{4}$},
$\alpha_{2} + \alpha_{4} + \alpha_{5}$,
$\alpha_{3} + \alpha_{4} + \alpha_{5}$,
\underline{$\alpha_{4} + \alpha_{5} + \alpha_{6}$},
$\alpha_{5} + \alpha_{6} + \alpha_{7}$,
\\[4pt]
\hline
&\\
\hline
\end{tabularx}
\end{landscape}
\begin{landscape}
\begin{tabularx}{\linewidth}{|>{\raggedright\let\newline\\\arraybackslash\hspace{0pt}}>{\hsize=1.35\hsize}X| >{\raggedright\let\newline\\\arraybackslash\hspace{0pt}}>{\hsize=0.65\hsize}X |}
\hline
&\\[-10pt]
$x_{47}=[[x_{7}, [x_{5}, x_{6}]], [[x_{3}, x_{4}], [x_{2}, [x_{1}, x_{3}]]]]$,
$x_{48}=[[[x_{1}, x_{3}], [x_{2}, x_{4}]], [[x_{2}, x_{3}], [x_{4}, x_{5}]]]$,
$\boxed{x_{49} }=[[[x_{1}, x_{3}], [x_{2}, x_{4}]], [[x_{3}, x_{4}], [x_{5}, x_{6}]]]$,
$x_{50}=[[[x_{1}, x_{3}], [x_{2}, x_{4}]], [[x_{4}, x_{5}], [x_{6}, x_{7}]]]$,
$x_{51}=[[[x_{1}, x_{3}], [x_{4}, x_{5}]], [[x_{2}, x_{4}], [x_{5}, x_{6}]]]$,
$x_{52}=[[[x_{2}, x_{3}], [x_{4}, x_{5}]], [[x_{4}, x_{5}], [x_{6}, x_{7}]]]$,
&
$2  \alpha_{3} + \alpha_{1} + \alpha_{2} + \alpha_{4} + \alpha_{5} + \alpha_{6} + \alpha_{7}$,
$2  \alpha_{2} + 2  \alpha_{3} + 2  \alpha_{4} + \alpha_{1} + \alpha_{5}$,
\underline{$2  \alpha_{3} + 2  \alpha_{4} + \alpha_{1} + \alpha_{2} + \alpha_{5} + \alpha_{6}$},
$2  \alpha_{4} + \alpha_{1} + \alpha_{2} + \alpha_{3} + \alpha_{5} + \alpha_{6} + \alpha_{7}$,
$2  \alpha_{4} + 2  \alpha_{5} + \alpha_{1} + \alpha_{2} + \alpha_{3} + \alpha_{6}$,
$2  \alpha_{4} + 2  \alpha_{5} + \alpha_{2} + \alpha_{3} + \alpha_{6} + \alpha_{7}$,
\\[4pt]
\hline
&\\
\hline
\end{tabularx}
\end{landscape}
\begin{landscape}
\begin{tabularx}{\linewidth}{|>{\raggedright\let\newline\\\arraybackslash\hspace{0pt}}>{\hsize=1.35\hsize}X| >{\raggedright\let\newline\\\arraybackslash\hspace{0pt}}>{\hsize=0.65\hsize}X |}
\hline
&\\[-10pt]
$x_{75}=[[[x_{2}, [x_{1}, x_{3}]], [[x_{3}, x_{4}], [x_{5}, x_{6}]]], [[x_{4}, [x_{2}, x_{3}]], [[x_{4}, x_{5}], [x_{6}, x_{7}]]]]$,
$x_{76}=[[[x_{5}, [x_{3}, x_{4}]], [[x_{1}, x_{3}], [x_{2}, x_{4}]]], [[x_{7}, [x_{5}, x_{6}]], [[x_{1}, x_{3}], [x_{2}, x_{4}]]]]$,
$\boxed{x_{77} }=[[[x_{6}, [x_{4}, x_{5}]], [[x_{2}, x_{3}], [x_{4}, x_{5}]]], [[x_{7}, [x_{5}, x_{6}]], [[x_{1}, x_{3}], [x_{2}, x_{4}]]]]$,
&
$2  \alpha_{2} + 2  \alpha_{5} + 2  \alpha_{6} + 3  \alpha_{3} + 3  \alpha_{4} + \alpha_{1} + \alpha_{7}$,
$2  \alpha_{1} + 2  \alpha_{2} + 2  \alpha_{5} + 3  \alpha_{3} + 3  \alpha_{4} + \alpha_{6} + \alpha_{7}$,
\underline{$2  \alpha_{2} + 2  \alpha_{3} + 2  \alpha_{6} + 3  \alpha_{4} + 3  \alpha_{5} + \alpha_{1} + \alpha_{7}$},
\\[4pt]
\hline
&\\[-10pt]
$x_{78}=[[[x_{6}, [x_{4}, x_{5}]], [[x_{2}, x_{3}], [x_{4}, x_{5}]]], [[x_{7}, [x_{5}, x_{6}]], [[x_{3}, x_{4}], [x_{2}, [x_{1}, x_{3}]]]]]$,
$x_{79}=[[[x_{6}, [x_{4}, x_{5}]], [[x_{2}, x_{3}], [x_{4}, x_{5}]]], [[[x_{1}, x_{3}], [x_{2}, x_{4}]], [[x_{4}, x_{5}], [x_{6}, x_{7}]]]]$,
$x_{80}=[[[x_{7}, [x_{5}, x_{6}]], [[x_{1}, x_{3}], [x_{2}, x_{4}]]], [[[x_{1}, x_{3}], [x_{2}, x_{4}]], [[x_{3}, x_{4}], [x_{5}, x_{6}]]]]$,
&
$2  \alpha_{2} + 2  \alpha_{6} + 3  \alpha_{3} + 3  \alpha_{4} + 3  \alpha_{5} + \alpha_{1} + \alpha_{7}$,
$2  \alpha_{2} + 2  \alpha_{3} + 2  \alpha_{6} + 3  \alpha_{5} + 4  \alpha_{4} + \alpha_{1} + \alpha_{7}$,
$2  \alpha_{1} + 2  \alpha_{2} + 2  \alpha_{5} + 2  \alpha_{6} + 3  \alpha_{3} + 3  \alpha_{4} + \alpha_{7}$,
\\[4pt]
\hline
&\\
\hline
\end{tabularx}
\end{landscape}

\subsection{$p=2$}{}~{}

\sssec{Notation \protect$\fA\oplus_c \fB$ needed to describe
$\fbgl(4; \alpha)$, $\fe(6, 6)$, $\fe(7,6)$, and $\fe(8,
1)$}\label{4cases} This notation describes the case where $\fA$ and
$\fB$ are nontrivial central extensions of the Lie algebras $\fa$
and $\fb$, respectively, and $\fA\oplus_c \fB$
--- a nontrivial central extension of $\fa\oplus \fb$ (or, perhaps,
a more complicated $\fa\subplus \fb$) with 1-dimensional center
spanned by $c$ --- is such that the restriction of the extension of
$\fa\oplus \fb$ to $\fa$ gives $\fA$ and that to $\fb$ gives $\fB$.

In these four cases, $\fg(A)_\ev$ is of the form
\[\fg(B)\oplus_c\fhei(2)\simeq \fg(B)\oplus\Span(X^+,X^-),\] where the
matrix $B$ is not invertible (so $\fg(B)$ has a grading element $d$
and a central element $c$), and where $X^+$, $X^-$ and $c$ span the
Heisenberg Lie algebra $\fhei(2)$. The brackets are:
\begin{equation}\label{strange}
\begin{array}{l}
{}[\fg^{(1)}(B), X^\pm]=0;\\
{}[d,X^\pm]=X^\pm; \qquad([d,X^\pm]=\alpha X^\pm\text{~for $\fbgl(4; \alpha)$})\\
{} [X^+,X^-]=c. \end{array}
\end{equation}

The odd part of $\fg(A)$ (at least in two of the four cases)
consists of two copies of the same $\fg(B)$-module $N$, the
operators $\ad_{X^\pm}$ permute these copies, and $\ad_{X^\pm}^2=0$,
so each of the operators maps one of the copies to the other, and
this other copy to zero.

\sssec{$\fbgl(3;\alpha)$, where $\alpha\neq 0, 1$; $\sdim=10/8|8$}\label{bgl3}   The roots of $\fg=\fbgl(3;\alpha)$ are the same as those of $\fosp(4|2;\alpha)$ (or, more correctly, $\fwk(3;\alpha)$), with the same division into even and odd ones; $\fg_\ev\simeq\fgl(3)\oplus \Kee Z$ and the $\fg_\ev$-module $\fg_\od$ is the sum of two irreducibles whose highest weight vectors are $x_7$ and $y_1$, where the roots corresponding to $x_i$ and $y_i$ are opposite.

We consider the   Cartan matrix and the Chevalley basis elements
\[

\quad
\begin{tabular}{|l|l|} 
\hline
&\\[-11pt]
 the root vectors&the roots\\[2pt]
\hline
&\\[-11pt]
\fbox{$x_1$}, \ \fbox{$x_2$},\ \fbox{$x_3$},&$\underline{\alpha_1}$, \ $\underline{\alpha_2}$,\ $\underline{\alpha_3}$\\[2pt]
$x_4=[x_1,\;x_2], \ x_5 =[x_2,\; x_3], $&$\alpha_1+\alpha_2, \ \alpha_2+\alpha_3,$\\[2pt]
$\fbox{$x_6$}=[x_3, [x_1,\;x_2]]$&$\underline{\alpha_1+\alpha_2+\alpha_3}$\\
$x_7=[[x_1, x_2], [x_2, x_3]]$&$\alpha_1+2\alpha_2+\alpha_3$\\[2pt]
\hline
\end{tabular}
\]

\sssec{$\fbgl(4;\alpha)$, where $\alpha\neq 0, 1$;
$\sdim =18|16$}\label{bgl4} The roots of $\fbgl(4;\alpha)$ are the same as those of $\fwk(4;\alpha)$,  but divided into even and odd ones: $\fbgl(4;\alpha)_\ev=\fgl(4)\oplus_c\fhei(2)$ and $\fbgl(4;\alpha)_\od=N\boxtimes\id$, where $N$ is an $8$-dimensional $\fgl(4)$ module, and $\id$ is the irreducible 2-dimensional $\fhei(2)$-module. We consider the   Cartan matrix and the
Chevalley basis elements
\[

\quad
\begin{tabular}{|l|l|}
  \hline
&\\[-11pt]
root vectors&roots\\[1pt]
  \hline
&\\[-11pt]
\fbox{$x_1$},\ \fbox{$x_2$} \ \fbox{$x_3$}, \ \fbox{$x_4$},&$\underline{\alpha_1}$, \ $\underline{\alpha_2}$, \ $\underline{\alpha_3}$, \ $\underline{\alpha_4}$\\[1pt]
$x_5=[x_1,\;x_2]$,\ $x_6=[x_1, x_3]$,
&$\alpha_1+\alpha_2, \ \alpha_1+\alpha_3$,
\\[1pt]
$x_7=[x_3,\;x_4]$,&$ \alpha_3+\alpha_4$\\[1pt]
$\fbox{$x_8$}=[x_3,\;[x_1,x_2]]$,
&$\underline{\alpha_1+\alpha_2+\alpha_3}$,
\\[1pt]
$\fbox{$x_9$}=[x_4,\;[x_1,\;x_3]]$,& $\underline{\alpha_1+\alpha_3+\alpha_4}$\\[1pt]
$x_{10}=[[x_1,x_2],\;[x_1,\;x_3]]$,
&$2\alpha_1+\alpha_2+\alpha_3$,
\\[1pt]
$ x_{11}=[[x_1,\;x_2],[x_3,\;x_4]]$& $ \alpha_1+\alpha_2+\alpha_3+\alpha_4 $\\[1pt]
\fbox{$x_{12}$}=$[[x_1,\;x_2],\;[x_4,[x_1,\;x_3]]],$&$\underline{2\alpha_1+\alpha_2+\alpha_3+\alpha_4}
$\\[1pt]
$x_{13}=[[x_3,\;[x_1,\;x_2]],[x_4,\;[x_1,\;x_3]]],$&$2\alpha_1+\alpha_2+2\alpha_3+\alpha_4$\\[1pt]
$\fbox{$x_{14}$}=[[x_4,\;[x_1,x_3]],\;[[x_1,\;x_2],\;[x_1\text{,$\,
 $}x_3]]]$&$\underline{3\alpha_1+\alpha_2+2\alpha_3+\alpha_4} $\\[2pt]
$x_{15}=[[[x_1,\;x_2],\;[x_1,\;x_3]],[[x_1,\;x_2],\;[x_3,x_4]]]$&$3\alpha_1+2\alpha_2+2\alpha_3+\alpha_4$\\[1pt]
  \hline
 \end{tabular}
\]

\ssec{$\fe(6,1)$ of $\sdim = 46|32$} We have
$\fg_\ev\simeq \fo\fc(2; 10)\oplus\Kee z$ and $\fg_\od$ is a
reducible module of the form $R(\pi_{4})\oplus R(\pi_{5})$ with the
two highest weight vectors \[x_{36}=[[[x_4,x_5],[x_6,[x_2,x_3]]],
[[x_3,[x_1,x_2]],[x_6,[x_3,x_4]]]]\] and $y_5$. We consider the Cartan matrix of $\fe(6)$ with parities of simple roots $111100$. The Chevalley basis elements are
$$
\small
\begin{tabularx}{\linewidth}{|>{\raggedright\let\newline\\\arraybackslash\hspace{0pt}}>{\hsize=1.25\hsize}X| >{\raggedright\let\newline\\\arraybackslash\hspace{0pt}}>{\hsize=0.75\hsize}X |}
\hline
&\\[-11pt]
the root vectors&the roots\\[1pt]
\hline
$\boxed{x_{1}}$,
$\boxed{x_{2}}$,
$\boxed{x_{3}}$,
$\boxed{x_{4}}$,
$x_{5}$,
$x_{6}$,
&
\underline{$\alpha_{1}$},
\underline{$\alpha_{2}$},
\underline{$\alpha_{3}$},
\underline{$\alpha_{4}$},
$\alpha_{5}$,
$\alpha_{6}$,
\\[1pt]
\hline
&\\[-11pt]
$x_{7}=[x_{1}, x_{2}]$,
$x_{8}=[x_{2}, x_{3}]$,
$x_{9}=[x_{3}, x_{4}]$,
$\boxed{x_{10} }=[x_{3}, x_{6}]$,
$\boxed{x_{11} }=[x_{4}, x_{5}]$,
&
$\alpha_{1} + \alpha_{2}$,
$\alpha_{2} + \alpha_{3}$,
$\alpha_{3} + \alpha_{4}$,
\underline{$\alpha_{3} + \alpha_{6}$},
\underline{$\alpha_{4} + \alpha_{5}$},
\\[1pt]
\hline
&\\[-11pt]
$\boxed{x_{12} }=[x_{3}, [x_{1}, x_{2}]]$,
$\boxed{x_{13} }=[x_{4}, [x_{2}, x_{3}]]$,
$x_{14}=[x_{5}, [x_{3}, x_{4}]]$,
$x_{15}=[x_{6}, [x_{2}, x_{3}]]$,\newline
$x_{16}=[x_{6}, [x_{3}, x_{4}]]$,
&
\underline{$\alpha_{1} + \alpha_{2} + \alpha_{3}$},
\underline{$\alpha_{2} + \alpha_{3} + \alpha_{4}$},
$\alpha_{3} + \alpha_{4} + \alpha_{5}$,
$\alpha_{2} + \alpha_{3} + \alpha_{6}$,
$\alpha_{3} + \alpha_{4} + \alpha_{6}$,
\\[1pt]
\hline
&\\
\hline
\end{tabularx}
$$

\sssec{$\fe(6,6)$ of $\sdim = 38|40$}  
In this case, $\fg(B)\simeq \fgl(6)$, see \eqref{4cases}. The $\fg_\ev$-module
$\fg_\od$ is irreducible with the highest weight vector
\[x_{35}=[[[x_3,x_6], [x_4,[x_2,x_3]]], [[x_4,x_5], [x_3,[x_1,x_2]]]]
\text{~~of weight $(0,0,1,0,0,1)$.}
\]
We consider the Cartan matrix of $\fe(6)$ with parities of simple roots $111111$. The Chevalley basis elements are
$$
\begin{tabularx}{\linewidth}{|>{\raggedright\let\newline\\\arraybackslash\hspace{0pt}}>{\hsize=1.25\hsize}X| >{\raggedright\let\newline\\\arraybackslash\hspace{0pt}}>{\hsize=0.75\hsize}X |}
\hline
&\\[-10pt]
the root vectors&the roots\\[4pt]
\hline
$\boxed{x_{1}}$,
$\boxed{x_{2}}$,
$\boxed{x_{3}}$,
$\boxed{x_{4}}$,
$\boxed{x_{5}}$,
$\boxed{x_{6}}$,
&
\underline{$\alpha_{1}$},
\underline{$\alpha_{2}$},
\underline{$\alpha_{3}$},
\underline{$\alpha_{4}$},
\underline{$\alpha_{5}$},
\underline{$\alpha_{6}$},
\\[4pt]
\hline
&\\[-10pt]
$x_{7}=[x_{1}, x_{2}]$,
$x_{8}=[x_{2}, x_{3}]$,
$x_{9}=[x_{3}, x_{4}]$,
$x_{10}=[x_{3}, x_{6}]$,
$x_{11}=[x_{4}, x_{5}]$,
&
$\alpha_{1} + \alpha_{2}$,
$\alpha_{2} + \alpha_{3}$,
$\alpha_{3} + \alpha_{4}$,
$\alpha_{3} + \alpha_{6}$,
$\alpha_{4} + \alpha_{5}$,
\\[4pt]
\hline
&\\[-10pt]
$\boxed{x_{12} }=[x_{3}, [x_{1}, x_{2}]]$,
$\boxed{x_{13} }=[x_{4}, [x_{2}, x_{3}]]$,
$\boxed{x_{14} }=[x_{5}, [x_{3}, x_{4}]]$,
$\boxed{x_{15} }=[x_{6}, [x_{2}, x_{3}]]$,
$\boxed{x_{16} }=[x_{6}, [x_{3}, x_{4}]]$,
&
\underline{$\alpha_{1} + \alpha_{2} + \alpha_{3}$},
\underline{$\alpha_{2} + \alpha_{3} + \alpha_{4}$},
\underline{$\alpha_{3} + \alpha_{4} + \alpha_{5}$},
\underline{$\alpha_{2} + \alpha_{3} + \alpha_{6}$},
\underline{$\alpha_{3} + \alpha_{4} + \alpha_{6}$},
\\[4pt]
\hline
&\\
\hline
\end{tabularx}
$$


\sssec{$\fe(7,1)$ of $\sdim=80/78|54$} We consider the Cartan matrix of $\fe(7)$ with the parities  of simple roots $1111001$. The Chevalley basis elements are 
$$
\footnotesize
\begin{tabularx}{\linewidth}{|>{\raggedright\let\newline\\\arraybackslash\hspace{0pt}}>{\hsize=1.23\hsize}X| >{\raggedright\let\newline\\\arraybackslash\hspace{0pt}}>{\hsize=0.77\hsize}X |}
\hline
&\\[-10pt]
the root vectors&the roots\\
\hline
$\boxed{x_{1}}$,
$\boxed{x_{2}}$,
$\boxed{x_{3}}$,
$\boxed{x_{4}}$,
$x_{5}$,
$x_{6}$,
$\boxed{x_{7}}$,
&
\underline{$\alpha_{1}$},
\underline{$\alpha_{2}$},
\underline{$\alpha_{3}$},
\underline{$\alpha_{4}$},
$\alpha_{5}$,
$\alpha_{6}$,
\underline{$\alpha_{7}$},
\\[4pt]
\hline
&\\[-10pt]
$x_{8}=[x_{1}, x_{2}]$,
$x_{9}=[x_{2}, x_{3}]$,
$x_{10}=[x_{3}, x_{4}]$,
$\boxed{x_{11} }=[x_{4}, x_{5}]$,
$x_{12}=[x_{4}, x_{7}]$,
$x_{13}=[x_{5}, x_{6}]$,
&
$\alpha_{1} + \alpha_{2}$,
$\alpha_{2} + \alpha_{3}$,
$\alpha_{3} + \alpha_{4}$,
\underline{$\alpha_{4} + \alpha_{5}$},
$\alpha_{4} + \alpha_{7}$,
$\alpha_{5} + \alpha_{6}$,
\\[4pt]
\hline
&\\[-10pt]
$\boxed{x_{14} }=[x_{3}, [x_{1}, x_{2}]]$,
$\boxed{x_{15} }=[x_{4}, [x_{2}, x_{3}]]$,\newline
$x_{16}=[x_{5}, [x_{3}, x_{4}]]$,
$\boxed{x_{17} }=[x_{6}, [x_{4}, x_{5}]]$,\newline
$\boxed{x_{18} }=[x_{7}, [x_{3}, x_{4}]]$,
$x_{19}=[x_{7}, [x_{4}, x_{5}]]$,
&
\underline{$\alpha_{1} + \alpha_{2} + \alpha_{3}$},
\underline{$\alpha_{2} + \alpha_{3} + \alpha_{4}$},
$\alpha_{3} + \alpha_{4} + \alpha_{5}$,
\underline{$\alpha_{4} + \alpha_{5} + \alpha_{6}$},
\underline{$\alpha_{3} + \alpha_{4} + \alpha_{7}$},
$\alpha_{4} + \alpha_{5} + \alpha_{7}$,
\\[4pt]
\hline
&\\
\hline
\end{tabularx}
$$

$$
\footnotesize
\begin{tabularx}{\linewidth}{|>{\raggedright\let\newline\\\arraybackslash\hspace{0pt}}>{\hsize=1.23\hsize}X| >{\raggedright\let\newline\\\arraybackslash\hspace{0pt}}>{\hsize=0.77\hsize}X |}
\hline
&\\[-10pt]
$\begin{array}{@{}r@{\mskip\medmuskip}c@{\mskip\medmuskip}l@{}}
x_{62}&=&[[[x_{7}, [x_{5}, [x_{3}, x_{4}]]], [[x_{1}, x_{2}], [x_{3}, x_{4}]]],\\
&& [[[x_{2}, x_{3}], [x_{4}, x_{5}]], [[x_{4}, x_{7}], [x_{5}, x_{6}]]]],
\end{array}
$
&
$2  \alpha_{2} + 2  \alpha_{7} + 3  \alpha_{3} + 3  \alpha_{5} + 4  \alpha_{4} + \alpha_{1} + \alpha_{6}$,
\\[4pt]
\hline
&\\[-10pt]
$\begin{array}{@{}r@{\mskip\medmuskip}c@{\mskip\medmuskip}l@{}}
x_{63}&=&[[[[x_{2}, x_{3}], [x_{4}, x_{7}]], [[x_{3}, x_{4}], [x_{5}, x_{6}]]],\\
&& [[[x_{4}, x_{7}], [x_{5}, x_{6}]], [[x_{4}, x_{5}], [x_{3}, [x_{1}, x_{2}]]]]]
\end{array}
$
&
$2  \alpha_{2} + 2  \alpha_{6} + 2  \alpha_{7} + 3  \alpha_{3} + 3  \alpha_{5} + 4  \alpha_{4} + \alpha_{1}$
\\[4pt]
\hline
\end{tabularx}
$$

\sssec{$\fe(7,6)$  of $\sdim=70/68|64$} We consider the Cartan matrix of $\fe(7)$ with the parities  of simple roots $0101010$ and the Chevalley basis elements
\begin{equation*}
\small
\begin{tabularx}{\linewidth}{|>{\raggedright\let\newline\\\arraybackslash\hspace{0pt}}>{\hsize=1.16\hsize}X| >{\raggedright\let\newline\\\arraybackslash\hspace{0pt}}>{\hsize=0.84\hsize}X |}
\hline
&\\[-10pt]
the root vectors&the roots\\[4pt]
\hline
&\\[-10pt]
$x_{1}$,
$\boxed{x_{2}}$,
$x_{3}$,
$\boxed{x_{4}}$,
$x_{5}$,
$\boxed{x_{6}}$,
$x_{7}$,
&
$\alpha_{1}$,
\underline{$\alpha_{2}$},
$\alpha_{3}$,
\underline{$\alpha_{4}$},
$\alpha_{5}$,
\underline{$\alpha_{6}$},
$\alpha_{7}$,
\\[4pt]
\hline
&\\[-10pt]
$\boxed{x_{8} }=[x_{1}, x_{2}]$,
$\boxed{x_{9} }=[x_{2}, x_{3}]$,
$\boxed{x_{10} }=[x_{3}, x_{4}]$,
$\boxed{x_{11} }=[x_{4}, x_{5}]$,
$\boxed{x_{12} }=[x_{4}, x_{7}]$,
$\boxed{x_{13} }=[x_{5}, x_{6}]$,
&
\underline{$\alpha_{1} + \alpha_{2}$},
\underline{$\alpha_{2} + \alpha_{3}$},
\underline{$\alpha_{3} + \alpha_{4}$},\newline
\underline{$\alpha_{4} + \alpha_{5}$},
\underline{$\alpha_{4} + \alpha_{7}$},
\underline{$\alpha_{5} + \alpha_{6}$},
\\[4pt]
\hline
&\\[-10pt]
$\boxed{x_{14} }=[x_{3}, [x_{1}, x_{2}]]$,
$x_{15}=[x_{4}, [x_{2}, x_{3}]]$,
$\boxed{x_{16} }=[x_{5}, [x_{3}, x_{4}]]$,
$x_{17}=[x_{6}, [x_{4}, x_{5}]]$,
$\boxed{x_{18} }=[x_{7}, [x_{3}, x_{4}]]$,
$\boxed{x_{19} }=[x_{7}, [x_{4}, x_{5}]]$,
&
\underline{$\alpha_{1} + \alpha_{2} + \alpha_{3}$},
$\alpha_{2} + \alpha_{3} + \alpha_{4}$,
\underline{$\alpha_{3} + \alpha_{4} + \alpha_{5}$},
$\alpha_{4} + \alpha_{5} + \alpha_{6}$,
\underline{$\alpha_{3} + \alpha_{4} + \alpha_{7}$},
\underline{$\alpha_{4} + \alpha_{5} + \alpha_{7}$},
\\
\hline
\end{tabularx}
\end{equation*}

\begin{equation*}
\small
\begin{tabularx}{\linewidth}{|>{\raggedright\let\newline\\\arraybackslash\hspace{0pt}}>{\hsize=1.16\hsize}X| >{\raggedright\let\newline\\\arraybackslash\hspace{0pt}}>{\hsize=0.84\hsize}X |}
\hline
&\\[-10pt]
$\begin{array}{@{}r@{\mskip\medmuskip}c@{\mskip\medmuskip}l@{}}
\boxed{x_{53} }&=&[[[x_{4}, x_{7}], [x_{5}, [x_{3}, x_{4}]]],\\
&& [[x_{3}, [x_{1}, x_{2}]], [x_{6}, [x_{4}, x_{5}]]]],
\end{array}
$
&
\underline{$2  \alpha_{3} + 2  \alpha_{5} + 3  \alpha_{4} + \alpha_{1} + \alpha_{2} + \alpha_{6} + \alpha_{7}$},\\
$\begin{array}{@{}r@{\mskip\medmuskip}c@{\mskip\medmuskip}l@{}}
\boxed{x_{54} }&=&[[[x_{5}, x_{6}], [x_{4}, [x_{2}, x_{3}]]],\\
&& [[x_{3}, [x_{1}, x_{2}]], [x_{7}, [x_{4}, x_{5}]]]],
\end{array}
$
&
\underline{$2  \alpha_{2} + 2  \alpha_{3} + 2  \alpha_{4} + 2  \alpha_{5} + \alpha_{1} + \alpha_{6} + \alpha_{7}$},\\
$\begin{array}{@{}r@{\mskip\medmuskip}c@{\mskip\medmuskip}l@{}}
\boxed{x_{55} }&=&[[[x_{5}, x_{6}], [x_{7}, [x_{3}, x_{4}]]],\\
&& [[x_{4}, [x_{2}, x_{3}]], [x_{7}, [x_{4}, x_{5}]]]],
\end{array}
$
&
\underline{$2  \alpha_{3} + 2  \alpha_{5} + 2  \alpha_{7} + 3  \alpha_{4} + \alpha_{2} + \alpha_{6}$},
\\[4pt]
\hline
&\\[-10pt]
$\begin{array}{@{}r@{\mskip\medmuskip}c@{\mskip\medmuskip}l@{}}
x_{56}&=&[[[x_{3}, [x_{1}, x_{2}]], [x_{6}, [x_{4}, x_{5}]]],\\
&& [[x_{4}, [x_{2}, x_{3}]], [x_{7}, [x_{4}, x_{5}]]]],
\end{array}
$
&
$2  \alpha_{2} + 2  \alpha_{3} + 2  \alpha_{5} + 3  \alpha_{4} + \alpha_{1} + \alpha_{6} + \alpha_{7}$,\\
$\begin{array}{@{}r@{\mskip\medmuskip}c@{\mskip\medmuskip}l@{}}
\boxed{x_{57} }&=&[[[x_{3}, [x_{1}, x_{2}]], [x_{7}, [x_{4}, x_{5}]]],\\
&& [[x_{6}, [x_{4}, x_{5}]], [x_{7}, [x_{3}, x_{4}]]]],
\end{array}
$
&
\underline{$2  \alpha_{3} + 2  \alpha_{5} + 2  \alpha_{7} + 3  \alpha_{4} + \alpha_{1} + \alpha_{2} + \alpha_{6}$},
\\[4pt]
\hline
&\\
\hline
\end{tabularx}
\end{equation*}

\sssec{$\fe(7,7)$  of $\sdim=64/62|70$} We consider the Cartan matrix of $\fe(7)$ with the parities  of simple roots $1111111$ and the Chevalley basis elements  
\begin{equation*}
\small
\begin{tabularx}{\linewidth}{|>{\raggedright\let\newline\\\arraybackslash\hspace{0pt}}>{\hsize=1.15\hsize}X| >{\raggedright\let\newline\\\arraybackslash\hspace{0pt}}>{\hsize=0.85\hsize}X |}
\hline
&\\[-10pt]
the root vectors&the roots\\[4pt]
\hline
&\\[-10pt]
$\boxed{x_{1}}$,
$\boxed{x_{2}}$,
$\boxed{x_{3}}$,
$\boxed{x_{4}}$,
$\boxed{x_{5}}$,
$\boxed{x_{6}}$,
$\boxed{x_{7}}$,
&
\underline{$\alpha_{1}$},
\underline{$\alpha_{2}$},
\underline{$\alpha_{3}$},
\underline{$\alpha_{4}$},
\underline{$\alpha_{5}$},
\underline{$\alpha_{6}$},
\underline{$\alpha_{7}$},
\\[4pt]
\hline
&\\[-10pt]
$x_{8}=[x_{1}, x_{2}]$,
$x_{9}=[x_{2}, x_{3}]$,
$x_{10}=[x_{3}, x_{4}]$,
$x_{11}=[x_{4}, x_{5}]$,
$x_{12}=[x_{4}, x_{7}]$,
$x_{13}=[x_{5}, x_{6}]$,
&
$\alpha_{1} + \alpha_{2}$,
$\alpha_{2} + \alpha_{3}$,
$\alpha_{3} + \alpha_{4}$,
$\alpha_{4} + \alpha_{5}$,
$\alpha_{4} + \alpha_{7}$,
$\alpha_{5} + \alpha_{6}$,
\\[4pt]
\hline
&\\[-10pt]
$\boxed{x_{14} }=[x_{3}, [x_{1}, x_{2}]]$,
$\boxed{x_{15} }=[x_{4}, [x_{2}, x_{3}]]$,
$\boxed{x_{16} }=[x_{5}, [x_{3}, x_{4}]]$,
$\boxed{x_{17} }=[x_{6}, [x_{4}, x_{5}]]$,
$\boxed{x_{18} }=[x_{7}, [x_{3}, x_{4}]]$,
$\boxed{x_{19} }=[x_{7}, [x_{4}, x_{5}]]$,
&
\underline{$\alpha_{1} + \alpha_{2} + \alpha_{3}$},
\underline{$\alpha_{2} + \alpha_{3} + \alpha_{4}$},
\underline{$\alpha_{3} + \alpha_{4} + \alpha_{5}$},
\underline{$\alpha_{4} + \alpha_{5} + \alpha_{6}$},
\underline{$\alpha_{3} + \alpha_{4} + \alpha_{7}$},
\underline{$\alpha_{4} + \alpha_{5} + \alpha_{7}$},
\\
\hline
\end{tabularx}
\end{equation*}

\begin{equation*}
\small
\begin{tabularx}{\linewidth}{|>{\raggedright\let\newline\\\arraybackslash\hspace{0pt}}>{\hsize=1.15\hsize}X| >{\raggedright\let\newline\\\arraybackslash\hspace{0pt}}>{\hsize=0.85\hsize}X |}
\hline
&\\[-10pt]
$x_{32}=[[x_{5}, x_{6}], [[x_{2}, x_{3}], [x_{4}, x_{7}]]]$,
$x_{33}=[[x_{3}, [x_{1}, x_{2}]], [x_{6}, [x_{4}, x_{5}]]]$,
$x_{34}=[[x_{3}, [x_{1}, x_{2}]], [x_{7}, [x_{4}, x_{5}]]]$,
$x_{35}=[[x_{4}, [x_{2}, x_{3}]], [x_{7}, [x_{4}, x_{5}]]]$,
$x_{36}=[[x_{6}, [x_{4}, x_{5}]], [x_{7}, [x_{3}, x_{4}]]]$,
&
$\alpha_{2} + \alpha_{3} + \alpha_{4} + \alpha_{5} + \alpha_{6} + \alpha_{7}$,
$\alpha_{1} + \alpha_{2} + \alpha_{3} + \alpha_{4} + \alpha_{5} + \alpha_{6}$,
$\alpha_{1} + \alpha_{2} + \alpha_{3} + \alpha_{4} + \alpha_{5} + \alpha_{7}$,
$2  \alpha_{4} + \alpha_{2} + \alpha_{3} + \alpha_{5} + \alpha_{7}$,
$2  \alpha_{4} + \alpha_{3} + \alpha_{5} + \alpha_{6} + \alpha_{7}$,
\\[4pt]
\hline
&\\[-10pt]
$\boxed{x_{37} }=[[x_{3}, [x_{1}, x_{2}]], [[x_{4}, x_{7}], [x_{5}, x_{6}]]]$,
$\boxed{x_{38} }=[[x_{6}, [x_{4}, x_{5}]], [[x_{2}, x_{3}], [x_{4}, x_{7}]]]$,
$\boxed{x_{39} }=[[x_{7}, [x_{3}, x_{4}]], [[x_{2}, x_{3}], [x_{4}, x_{5}]]]$,
$\boxed{x_{40} }=[[x_{7}, [x_{4}, x_{5}]], [[x_{1}, x_{2}], [x_{3}, x_{4}]]]$,
$\boxed{x_{41} }=[[x_{7}, [x_{4}, x_{5}]], [[x_{3}, x_{4}], [x_{5}, x_{6}]]]$,
&
\underline{$\alpha_{1} + \alpha_{2} + \alpha_{3} + \alpha_{4} + \alpha_{5} + \alpha_{6} + \alpha_{7}$},
\underline{$2  \alpha_{4} + \alpha_{2} + \alpha_{3} + \alpha_{5} + \alpha_{6} + \alpha_{7}$},
\underline{$2  \alpha_{3} + 2  \alpha_{4} + \alpha_{2} + \alpha_{5} + \alpha_{7}$},
\underline{$2  \alpha_{4} + \alpha_{1} + \alpha_{2} + \alpha_{3} + \alpha_{5} + \alpha_{7}$},
\underline{$2  \alpha_{4} + 2  \alpha_{5} + \alpha_{3} + \alpha_{6} + \alpha_{7}$},
\\[4pt]
\hline
&\\
\hline
\end{tabularx}
\end{equation*}

\sssec{$\fe(8,1)$  of $\sdim=136|112$} We have (cf. Sect. \ref{4cases}) $\fg(B)\simeq \fe(7)$.  We consider the Cartan Matrix of $\mathfrak{e}(8)$ with the parities of simple roots $11001111$ and the Chevalley basis elements

\begin{landscape}
\begin{tabularx}{.97\linewidth}{|>{\raggedright\let\newline\\\arraybackslash\hspace{0pt}}>{\hsize=1.2\hsize}X| >{\raggedright\let\newline\\\arraybackslash\hspace{0pt}}>{\hsize=0.8\hsize}X |}
\hline
&\\[-10pt]
the root vectors&the roots\\[4pt]
\hline
$\boxed{x_{1}}$,
$\boxed{x_{2}}$,
$x_{3}$,
$x_{4}$,
$\boxed{x_{5}}$,
$\boxed{x_{6}}$,
$\boxed{x_{7}}$,
$\boxed{x_{8}}$,
&
\underline{$\alpha_{1}$},
\underline{$\alpha_{2}$},
$\alpha_{3}$,
$\alpha_{4}$,
\underline{$\alpha_{5}$},
\underline{$\alpha_{6}$},
\underline{$\alpha_{7}$},
\underline{$\alpha_{8}$},
\\[4pt]
\hline
&\\[-10pt]
$x_{9}=[x_{1}, x_{2}]$,
$\boxed{x_{10} }=[x_{2}, x_{3}]$,
$x_{11}=[x_{3}, x_{4}]$,
$\boxed{x_{12} }=[x_{4}, x_{5}]$,
$x_{13}=[x_{5}, x_{6}]$,
$x_{14}=[x_{5}, x_{8}]$,
$x_{15}=[x_{6}, x_{7}]$,
&
$\alpha_{1} + \alpha_{2}$,
\underline{$\alpha_{2} + \alpha_{3}$},
$\alpha_{3} + \alpha_{4}$,
\underline{$\alpha_{4} + \alpha_{5}$},
$\alpha_{5} + \alpha_{6}$,
$\alpha_{5} + \alpha_{8}$,
$\alpha_{6} + \alpha_{7}$,
\\[4pt]
\hline
&\\[-10pt]
$x_{16}=[x_{3}, [x_{1}, x_{2}]]$,
$\boxed{x_{17} }=[x_{4}, [x_{2}, x_{3}]]$,
$\boxed{x_{18} }=[x_{5}, [x_{3}, x_{4}]]$,
$x_{19}=[x_{6}, [x_{4}, x_{5}]]$,
$\boxed{x_{20} }=[x_{7}, [x_{5}, x_{6}]]$,
$x_{21}=[x_{8}, [x_{4}, x_{5}]]$,
$\boxed{x_{22} }=[x_{8}, [x_{5}, x_{6}]]$,
&
$\alpha_{1} + \alpha_{2} + \alpha_{3}$,
\underline{$\alpha_{2} + \alpha_{3} + \alpha_{4}$},
\underline{$\alpha_{3} + \alpha_{4} + \alpha_{5}$},
$\alpha_{4} + \alpha_{5} + \alpha_{6}$,
\underline{$\alpha_{5} + \alpha_{6} + \alpha_{7}$},
$\alpha_{4} + \alpha_{5} + \alpha_{8}$,
\underline{$\alpha_{5} + \alpha_{6} + \alpha_{8}$},
\\[4pt]
\hline
&\\
\hline
\end{tabularx}
\end{landscape}

\begin{landscape}
\begin{tabularx}{.97\linewidth}{|>{\raggedright\let\newline\\\arraybackslash\hspace{0pt}}>{\hsize=1.2\hsize}X| >{\raggedright\let\newline\\\arraybackslash\hspace{0pt}}>{\hsize=0.8\hsize}X |}
\hline
&\\[-10pt]
$x_{51}=[[x_{8}, [x_{5}, x_{6}]], [[x_{4}, x_{5}], [x_{3}, [x_{1}, x_{2}]]]]$,
$\boxed{x_{52} }=[[x_{8}, [x_{6}, [x_{4}, x_{5}]]], [[x_{2}, x_{3}], [x_{4}, x_{5}]]]$,
$x_{53}=[[[x_{1}, x_{2}], [x_{3}, x_{4}]], [[x_{5}, x_{8}], [x_{6}, x_{7}]]]$,
$x_{54}=[[[x_{2}, x_{3}], [x_{4}, x_{5}]], [[x_{5}, x_{8}], [x_{6}, x_{7}]]]$,
$x_{55}=[[[x_{3}, x_{4}], [x_{5}, x_{6}]], [[x_{5}, x_{8}], [x_{6}, x_{7}]]]$,
$\boxed{x_{56} }=[[[x_{3}, x_{4}], [x_{5}, x_{8}]], [[x_{4}, x_{5}], [x_{6}, x_{7}]]]$,
&
$2  \alpha_{5} + \alpha_{1} + \alpha_{2} + \alpha_{3} + \alpha_{4} + \alpha_{6} + \alpha_{8}$,
\underline{$2  \alpha_{4} + 2  \alpha_{5} + \alpha_{2} + \alpha_{3} + \alpha_{6} + \alpha_{8}$},
$\alpha_{1} + \alpha_{2} + \alpha_{3} + \alpha_{4} + \alpha_{5} + \alpha_{6} + \alpha_{7} + \alpha_{8}$,
$2  \alpha_{5} + \alpha_{2} + \alpha_{3} + \alpha_{4} + \alpha_{6} + \alpha_{7} + \alpha_{8}$,
$2  \alpha_{5} + 2  \alpha_{6} + \alpha_{3} + \alpha_{4} + \alpha_{7} + \alpha_{8}$,
\underline{$2  \alpha_{4} + 2  \alpha_{5} + \alpha_{3} + \alpha_{6} + \alpha_{7} + \alpha_{8}$},
\\[4pt]
\hline
&\\
\hline
\end{tabularx}
\end{landscape}

\begin{landscape}
\begin{tabularx}{.97\linewidth}{|>{\raggedright\let\newline\\\arraybackslash\hspace{0pt}}>{\hsize=1.2\hsize}X| >{\raggedright\let\newline\\\arraybackslash\hspace{0pt}}>{\hsize=0.8\hsize}X |}
\hline
&\\[-10pt]
$x_{80}=[[[x_{4}, [x_{2}, x_{3}]], [x_{7}, [x_{5}, x_{6}]]], [[x_{8}, [x_{4}, x_{5}]], [[x_{3}, x_{4}], [x_{5}, x_{6}]]]]$,
$\boxed{x_{81} }=[[[x_{4}, [x_{2}, x_{3}]], [x_{8}, [x_{5}, x_{6}]]], [[x_{7}, [x_{5}, x_{6}]], [[x_{1}, x_{2}], [x_{3}, x_{4}]]]]$,
$\boxed{x_{82} }=[[[x_{5}, [x_{3}, x_{4}]], [x_{8}, [x_{5}, x_{6}]]], [[x_{4}, [x_{2}, x_{3}]], [[x_{5}, x_{8}], [x_{6}, x_{7}]]]]$,
$\boxed{x_{83} }=[[[x_{5}, [x_{3}, x_{4}]], [x_{8}, [x_{5}, x_{6}]]], [[x_{7}, [x_{5}, x_{6}]], [[x_{1}, x_{2}], [x_{3}, x_{4}]]]]$,
$x_{84}=[[[x_{7}, [x_{5}, x_{6}]], [x_{8}, [x_{4}, x_{5}]]], [[x_{8}, [x_{5}, x_{6}]], [[x_{1}, x_{2}], [x_{3}, x_{4}]]]]$,
&
$2  \alpha_{3} + 2  \alpha_{6} + 3  \alpha_{4} + 3  \alpha_{5} + \alpha_{2} + \alpha_{7} + \alpha_{8}$,
\underline{$2  \alpha_{2} + 2  \alpha_{3} + 2  \alpha_{4} + 2  \alpha_{5} + 2  \alpha_{6} + \alpha_{1} + \alpha_{7} + \alpha_{8}$},
\underline{$2  \alpha_{3} + 2  \alpha_{4} + 2  \alpha_{6} + 2  \alpha_{8} + 3  \alpha_{5} + \alpha_{2} + \alpha_{7}$},
\underline{$2  \alpha_{3} + 2  \alpha_{4} + 2  \alpha_{6} + 3  \alpha_{5} + \alpha_{1} + \alpha_{2} + \alpha_{7} + \alpha_{8}$},
$2  \alpha_{4} + 2  \alpha_{6} + 2  \alpha_{8} + 3  \alpha_{5} + \alpha_{1} + \alpha_{2} + \alpha_{3} + \alpha_{7}$,
\\[4pt]
\hline
&\\
\hline
\end{tabularx}
\end{landscape}

\begin{landscape}
\begin{tabularx}{.97\linewidth}{|>{\raggedright\let\newline\\\arraybackslash\hspace{0pt}}>{\hsize=1.2\hsize}X| >{\raggedright\let\newline\\\arraybackslash\hspace{0pt}}>{\hsize=0.8\hsize}X |}
\hline
&\\[-10pt]
$\begin{array}{@{}r@{\mskip\medmuskip}c@{\mskip\medmuskip}l@{}}
\boxed{x_{101} }&=&[[[[x_{1}, x_{2}], [x_{3}, x_{4}]], [[x_{5}, x_{8}], [x_{6}, [x_{4}, x_{5}]]]],\\
&& [[[x_{5}, x_{8}], [x_{6}, x_{7}]], [[x_{5}, x_{6}], [x_{4}, [x_{2}, x_{3}]]]]],
\end{array}
$
&
\underline{$2  \alpha_{2} + 2  \alpha_{3} + 2  \alpha_{8} + 3  \alpha_{4} + 3  \alpha_{6} + 4  \alpha_{5} + \alpha_{1} + \alpha_{7}$},\\
$\begin{array}{@{}r@{\mskip\medmuskip}c@{\mskip\medmuskip}l@{}}
x_{102}&=&[[[[x_{3}, x_{4}], [x_{5}, x_{8}]], [[x_{5}, x_{6}], [x_{4}, [x_{2}, x_{3}]]]],\\
&& [[[x_{5}, x_{8}], [x_{6}, x_{7}]], [[x_{4}, x_{5}], [x_{3}, [x_{1}, x_{2}]]]]],
\end{array}
$
&
$2  \alpha_{2} + 2  \alpha_{6} + 2  \alpha_{8} + 3  \alpha_{3} + 3  \alpha_{4} + 4  \alpha_{5} + \alpha_{1} + \alpha_{7}$,\\
$\begin{array}{@{}r@{\mskip\medmuskip}c@{\mskip\medmuskip}l@{}}
\boxed{x_{103} }&=&[[[[x_{4}, x_{5}], [x_{6}, x_{7}]], [[x_{3}, x_{4}], [x_{8}, [x_{5}, x_{6}]]]],\\
&& [[[x_{5}, x_{8}], [x_{6}, x_{7}]], [[x_{4}, x_{5}], [x_{3}, [x_{1}, x_{2}]]]]],
\end{array}
$
&
\underline{$2  \alpha_{3} + 2  \alpha_{7} + 2  \alpha_{8} + 3  \alpha_{4} + 3  \alpha_{6} + 4  \alpha_{5} + \alpha_{1} + \alpha_{2}$},
\\[4pt]
\hline
&\\[-10pt]
$\begin{array}{@{}r@{\mskip\medmuskip}c@{\mskip\medmuskip}l@{}}
x_{104}&=&[[[[x_{4}, x_{5}], [x_{6}, x_{7}]], [[x_{5}, x_{8}], [x_{4}, [x_{2}, x_{3}]]]],\\
&& [[[x_{3}, x_{4}], [x_{8}, [x_{5}, x_{6}]]], [[x_{4}, x_{5}], [x_{3}, [x_{1}, x_{2}]]]]],
\end{array}
$
&
$2  \alpha_{2} + 2  \alpha_{6} + 2  \alpha_{8} + 3  \alpha_{3} + 4  \alpha_{4} + 4  \alpha_{5} + \alpha_{1} + \alpha_{7}$,\\
\hline
\end{tabularx}
\end{landscape}

\begin{landscape}
\begin{tabularx}{.97\linewidth}{|>{\raggedright\let\newline\\\arraybackslash\hspace{0pt}}>{\hsize=1.2\hsize}X| >{\raggedright\let\newline\\\arraybackslash\hspace{0pt}}>{\hsize=0.8\hsize}X |}
\hline
&\\[-10pt]
$\begin{array}{@{}r@{\mskip\medmuskip}c@{\mskip\medmuskip}l@{}}
x_{116}&=&[[[[x_{4}, [x_{2}, x_{3}]], [x_{8}, [x_{5}, x_{6}]]], [[x_{7}, [x_{5}, x_{6}]], [x_{8}, [x_{4}, x_{5}]]]],\\
&& [[[x_{5}, [x_{3}, x_{4}]], [x_{8}, [x_{5}, x_{6}]]], [[x_{7}, [x_{5}, x_{6}]], [[x_{1}, x_{2}], [x_{3}, x_{4}]]]]],
\end{array}
$
&
$2  \alpha_{2} + 2  \alpha_{7} + 3  \alpha_{3} + 3  \alpha_{8} + 4  \alpha_{4} + 4  \alpha_{6} + 6  \alpha_{5} + \alpha_{1}$,
\\[4pt]
\hline
&\\[-10pt]
$\begin{array}{@{}r@{\mskip\medmuskip}c@{\mskip\medmuskip}l@{}}
x_{117}&=&[[[[x_{4}, [x_{2}, x_{3}]], [x_{7}, [x_{5}, x_{6}]]], [[x_{8}, [x_{4}, x_{5}]], [[x_{3}, x_{4}], [x_{5}, x_{6}]]]],\\
&& [[[x_{7}, [x_{5}, x_{6}]], [x_{8}, [x_{4}, x_{5}]]], [[x_{8}, [x_{5}, x_{6}]], [[x_{1}, x_{2}], [x_{3}, x_{4}]]]]],
\end{array}
$
&
$2  \alpha_{2} + 2  \alpha_{7} + 3  \alpha_{3} + 3  \alpha_{8} + 4  \alpha_{6} + 5  \alpha_{4} + 6  \alpha_{5} + \alpha_{1}$,
\\[4pt]
\hline
&\\[-10pt]
$\begin{array}{@{}r@{\mskip\medmuskip}c@{\mskip\medmuskip}l@{}}
x_{118}&=&[[[[x_{5}, [x_{3}, x_{4}]], [x_{8}, [x_{5}, x_{6}]]], [[x_{7}, [x_{5}, x_{6}]], [[x_{1}, x_{2}], [x_{3}, x_{4}]]]],\\
&& [[[x_{4}, [x_{2}, x_{3}]], [[x_{5}, x_{8}], [x_{6}, x_{7}]]], [[x_{8}, [x_{4}, x_{5}]], [[x_{3}, x_{4}], [x_{5}, x_{6}]]]]],
\end{array}
$
&
$2  \alpha_{2} + 2  \alpha_{7} + 3  \alpha_{8} + 4  \alpha_{3} + 4  \alpha_{6} + 5  \alpha_{4} + 6  \alpha_{5} + \alpha_{1}$,
\\[4pt]
\hline
&\\
\hline
\end{tabularx}
\end{landscape}

\sssec{$\fe(8,8)$ of $\sdim=120|128$} In the $\Zee$-grading
with the 1st CM with $\deg e_8^\pm=\pm 1$ and $\deg e_i^\pm=0$ for
$i\neq 8$, we have $\fg_0=\fgl(8)=\fgl(V)$. There are different
isomorphisms between $\fg_0$ and $\fgl(8)$; using the one where
$h_i=E_{i,i}+E_{i+1,i+1}$ for all $i=1,\dots, 7$, and
$h_8=E_{6,6}+E_{7,7}+E_{8,8}$, we see that, as modules over
$\fgl(V)$,
\[
\begin{array}{lll}
\fg_1=\bigwedge^5V^*;&\fg_2=\bigwedge^6 V;&\fg_3=V;\\
\fg_{-1}=\bigwedge^5V;&\fg_{-2}=\bigwedge^6 V^*;&\fg_{-3}=V^*.
\end{array}
\]
We can also set $h_8=E_{1,1}+E_{2,2}+E_{3,3}+E_{4,4}+E_{5,5}$. Then
we get
\[
\begin{array}{lll}
\fg_1=\bigwedge^3V;&\fg_2=\bigwedge^6 V;&\fg_3=\bigwedge^7V^*;\\
\fg_{-1}=\bigwedge^3V^*;&\fg_{-2}=\bigwedge^6
V^*;&\fg_{-3}=\bigwedge^7V.
\end{array}
\]

The algebra $\fg_\ev$ is isomorphic to
$\fo_\Pi^{(2)}(16)\subplus\Kee d$, where
$d=E_{6,6}+\dots+E_{13,13}$, and $\fg_\od$ is an irreducible
$\fg_\ev$-module with the highest weight element
$x_{120}$ of weight $(1,0,\dots,0)$ with respect to $h_1,\dots,h_8$;
$\fg_\od$ also possesses a lowest weight vector.

\begin{tabularx}{\linewidth}{|>{\raggedright\let\newline\\\arraybackslash\hspace{0pt}}>{\hsize=1.3\hsize}X| >{\raggedright\let\newline\\\arraybackslash\hspace{0pt}}>{\hsize=0.7\hsize}X |}
\hline
&\\[-10pt]
the root vectors&the roots\\[4pt]
\hline
$\boxed{x_{1}}$,
$\boxed{x_{2}}$,
$\boxed{x_{3}}$,
$\boxed{x_{4}}$,
$\boxed{x_{5}}$,
$\boxed{x_{6}}$,
$\boxed{x_{7}}$,
$\boxed{x_{8}}$,
&
\underline{$\alpha_{1}$},
\underline{$\alpha_{2}$},
\underline{$\alpha_{3}$},
\underline{$\alpha_{4}$},
\underline{$\alpha_{5}$},
\underline{$\alpha_{6}$},
\underline{$\alpha_{7}$},
\underline{$\alpha_{8}$},
\\[4pt]
\hline
&\\[-10pt]
$x_{9}=[x_{1}, x_{2}]$,
$x_{10}=[x_{2}, x_{3}]$,
$x_{11}=[x_{3}, x_{4}]$,
$x_{12}=[x_{4}, x_{5}]$,
$x_{13}=[x_{5}, x_{6}]$,
$x_{14}=[x_{5}, x_{8}]$,
$x_{15}=[x_{6}, x_{7}]$,
&
$\alpha_{1} + \alpha_{2}$,
$\alpha_{2} + \alpha_{3}$,
$\alpha_{3} + \alpha_{4}$,
$\alpha_{4} + \alpha_{5}$,
$\alpha_{5} + \alpha_{6}$,
$\alpha_{5} + \alpha_{8}$,
$\alpha_{6} + \alpha_{7}$,
\\[4pt]
\hline
&\\[-10pt]
$\boxed{x_{16} }=[x_{3}, [x_{1}, x_{2}]]$,
$\boxed{x_{17} }=[x_{4}, [x_{2}, x_{3}]]$,
$\boxed{x_{18} }=[x_{5}, [x_{3}, x_{4}]]$,
$\boxed{x_{19} }=[x_{6}, [x_{4}, x_{5}]]$,
$\boxed{x_{20} }=[x_{7}, [x_{5}, x_{6}]]$,
$\boxed{x_{21} }=[x_{8}, [x_{4}, x_{5}]]$,
$\boxed{x_{22} }=[x_{8}, [x_{5}, x_{6}]]$,
&
\underline{$\alpha_{1} + \alpha_{2} + \alpha_{3}$},
\underline{$\alpha_{2} + \alpha_{3} + \alpha_{4}$},
\underline{$\alpha_{3} + \alpha_{4} + \alpha_{5}$},
\underline{$\alpha_{4} + \alpha_{5} + \alpha_{6}$},
\underline{$\alpha_{5} + \alpha_{6} + \alpha_{7}$},
\underline{$\alpha_{4} + \alpha_{5} + \alpha_{8}$},
\underline{$\alpha_{5} + \alpha_{6} + \alpha_{8}$},
\\[4pt]
\hline
&\\
\hline
\end{tabularx}

\begin{landscape}
\begin{tabularx}{.97\linewidth}{|>{\raggedright\let\newline\\\arraybackslash\hspace{0pt}}>{\hsize=1.3\hsize}X| >{\raggedright\let\newline\\\arraybackslash\hspace{0pt}}>{\hsize=0.7\hsize}X |}
\hline
&\\[-10pt]
$\boxed{x_{44} }=[[x_{4}, [x_{2}, x_{3}]], [[x_{5}, x_{8}], [x_{6}, x_{7}]]]$,\newline
$\boxed{x_{45} }=[[x_{7}, [x_{5}, x_{6}]], [[x_{1}, x_{2}], [x_{3}, x_{4}]]]$,\newline
$\boxed{x_{46} }=[[x_{7}, [x_{5}, x_{6}]], [[x_{3}, x_{4}], [x_{5}, x_{8}]]]$,\newline
$\boxed{x_{47} }=[[x_{8}, [x_{4}, x_{5}]], [[x_{3}, x_{4}], [x_{5}, x_{6}]]]$,\newline
$\boxed{x_{48} }=[[x_{8}, [x_{5}, x_{6}]], [[x_{1}, x_{2}], [x_{3}, x_{4}]]]$,\newline
$\boxed{x_{49} }=[[x_{8}, [x_{5}, x_{6}]], [[x_{2}, x_{3}], [x_{4}, x_{5}]]]$,\newline
$\boxed{x_{50} }=[[x_{8}, [x_{5}, x_{6}]], [[x_{4}, x_{5}], [x_{6}, x_{7}]]]$,
&
\underline{$\alpha_{2} + \alpha_{3} + \alpha_{4} + \alpha_{5} + \alpha_{6} + \alpha_{7} + \alpha_{8}$},
\underline{$\alpha_{1} + \alpha_{2} + \alpha_{3} + \alpha_{4} + \alpha_{5} + \alpha_{6} + \alpha_{7}$},
\underline{$2  \alpha_{5} + \alpha_{3} + \alpha_{4} + \alpha_{6} + \alpha_{7} + \alpha_{8}$},
\underline{$2  \alpha_{4} + 2  \alpha_{5} + \alpha_{3} + \alpha_{6} + \alpha_{8}$},
\underline{$\alpha_{1} + \alpha_{2} + \alpha_{3} + \alpha_{4} + \alpha_{5} + \alpha_{6} + \alpha_{8}$},
\underline{$2  \alpha_{5} + \alpha_{2} + \alpha_{3} + \alpha_{4} + \alpha_{6} + \alpha_{8}$},
\underline{$2  \alpha_{5} + 2  \alpha_{6} + \alpha_{4} + \alpha_{7} + \alpha_{8}$},
\\[4pt]
\hline
&\\
\hline
\end{tabularx}%
\end{landscape}

\section{Root systems of  Lie algebras of the form $\fg(A)$ with indecomposable $A$}

\ssec{$\fwk(3; a)$ and $\fwk(4; a)$} These Lie algebras are desuperizations of $\fbgl(3; a)$ and $\fbgl(4; a)$, respectively, so they have the same root systems, see~Subsections~\ref{bgl3} and~\ref{bgl4}.

\ssec{$\textbf{F}(\fo\fo(1|2n))$, where $\textbf{F}$ is the desuperization functor} (In \cite{WK}, this simple Lie algebra is denoted $\Delta_n$.) Its root system is the same as that of $\fo(2n+1)$,  see \cite{BGLLS}.

\ssec{$\fbr(2; \eps)$, $p=3$} Observe that it has the same root system as $\fo(5)$.

\ssec{$\fbr(3)$}\label{br3} We consider the
Cartan matrix (Skryabin was the first to describe the Cartan matrices of $\fbr(3)$, see \cite{SkB}) and the Chevalley basis elements
\[
\begin{array}{ll}
\begin{pmatrix}
2&-1&0\\
-1&2&-1\\
0&-1&\ev\end{pmatrix}&
\begin{tabular}{|l|l|} 
\hline
 the root vectors&the roots\\
\hline
$x_1,\; x_2,\; x_3$&$\alpha_1,\; \alpha_2,\; \alpha_3$\\
$x_4=[x_1,x_2],\quad  x_5=[x_2,x_3]$&$\alpha_1+\alpha_2,\; \alpha_2+\alpha_3$\\
$x_6=[x_3,[x_1, x_2]],\quad x_7=[x_3,[x_2,x_3]]$&$\alpha_1+\alpha_2+\alpha_3,\; \alpha_2+2\alpha_3$\\
$x_8=[x_3,[x_3,[x_1,x_2]]]$&$\alpha_1+\alpha_2+2\alpha_3$\\
$x_9=[[x_2,x_3], [x_3,[x_1,x_2]]]$&$\alpha_1+2\alpha_2+2\alpha_3$\\
$x_{10}=[[x_3,[x_1,x_2]], [x_3,[x_2,x_3]]]$&$\alpha_1+2\alpha_2+3\alpha_3$\\
$x_{11}=[[x_3,[x_2,x_3]], [x_3,[x_3,[x_1,x_2]]]]$&$\alpha_1+2\alpha_2+4\alpha_3$\\
$x_{12}=[[x_3,[x_2,x_3]], [[x_2,x_3],
[x_3,[x_1,x_2]]]]$&$\alpha_1+3\alpha_2+4\alpha_3$\\
$x_{13}=[[x_3,[x_3,[x_1,x_2]]], [[x_2,x_3],
[x_3,[x_1,x_2]]]]$&$2\alpha_1+3\alpha_2+4\alpha_3$\\
\hline
\end{tabular}
\end{array}
\]


\end{document}